\documentclass[11pt]{amsart}
\usepackage{amscd,amssymb,hyperref}
\input xy
\xyoption{all}

\newtheorem{theorem}{Theorem}[section]

\newtheorem{conjecture}[theorem]{Conjecture}

\renewcommand{\leq}{\leqslant}
\renewcommand{\geq}{\geqslant}

\newtheorem{question}[theorem]{Question}
\theoremstyle{definition}
\newtheorem{definition}[theorem]{Definition}
\newtheorem{example}[theorem]{Example}

\theoremstyle{definition}

\numberwithin{equation}{section}

\newcommand{\ve}{\varepsilon}
\newcommand{\vp}{\varphi}
\newcommand{\mn}{\sqrt{-1}}
\newcommand{\ov}[1]{\overline{#1}}
\newcommand{\de}{\partial}
\newcommand{\db}{\overline{\partial}}
\newcommand{\ddbar}{\sqrt{-1} \partial \overline{\partial}}
\newcommand{\ti}[1]{\tilde{#1}}
\DeclareMathOperator{\Vol}{Vol}

\DeclareMathOperator{\Ric}{Ric} 
\newcommand{\diam}{\mathrm{diam}}
\newcommand{\Null}{\mathrm{Null}}

\numberwithin{equation}{section} \numberwithin{figure}{section}

\author{Valentino Tosatti}
\address{Department of Mathematics, Northwestern University, 2033 Sheridan Road, Evanston, IL 60208}
\email{tosatti@math.northwestern.edu}
\title{Collapsing Calabi-Yau manifolds}
\dedicatory{Dedicated to Professor Shing-Tung Yau on the occasion of his 70th birthday.}

\begin{document}

\begin{abstract}We survey some recent developments on the problem of understanding degenerations of Calabi-Yau manifolds equipped with their Ricci-flat K\"ahler metrics, with an emphasis on the case when the metrics are volume collapsing.
\end{abstract}

\maketitle

\section{Introduction}

Calabi-Yau manifolds form an important class of compact complex manifolds that enjoys remarkable geometric properties, and have been extensively studied in many fields of mathematics (as well as theoretical physics). The main feature of these manifolds is that they carry K\"ahler metrics with everywhere vanishing Ricci curvature, thanks to Yau's solution \cite{Ya} of the Calabi Conjecture \cite{Cal}. These Ricci-flat K\"ahler metrics are not flat (unless we are on a torus or a finite quotient), and do not have an explicit description; rather, they are constructed by solving a fully nonlinear PDE of complex Monge-Amp\`ere type.

In recent years, the study of families of such Ricci-flat Calabi-Yau manifolds has appeared naturally in many contexts, including in the study of moduli spaces and in mirror symmetry, and it is especially interesting to understand the possible ways in which these metrics can degenerate.

In this survey we will give an overview of some aspects of this problem, with a focus on the topics that are closest to the author's interests. Greater emphasis is placed on the case when the Ricci-flat manifolds are volume collapsing, which is by far the hardest case to analyze, hence the title of this survey. We will also describe a number of open problems and conjectures, some well-known and some less so, in the hope to stimulate further progress.

The reader is also encouraged to consult the author's earlier surveys \cite{To2,To3}, as well as the more recent surveys by Zhang \cite{Zh2} and Sun \cite{Su} on this circle of ideas. The relations with holomorphic dynamics that we will discuss in section \ref{below} is also described in the recent expository papers \cite{Fil} and \cite{ToS}.\\

{\bf Acknowledgements. }The author would like to thank Professor Shing-Tung Yau for many useful conversations on the topics of this article over the years, and for his advice and support. This survey is an expanded version of the author's talk at the conference ``Differential Geometry, Calabi-Yau Theory and General Relativity'' that took place at Harvard University on May 2--5, 2019, in celebration of Professor Shing-Tung Yau's 70th birthday. He is grateful to the organizers of the conference, and to the editors of Surveys in Differential Geometry for the invitation to submit this article.
The author was partially supported by NSF grant DMS-1903147. This article was written during the author's stay at the Center for Mathematical Sciences and Applications at Harvard University, which he would like to thank for the hospitality.

\section{Degenerations of Calabi-Yau manifolds}
\subsection{Calabi-Yau manifolds}
The main object of study in this article are compact Calabi-Yau manifolds. There are several slightly different definitions of Calabi-Yau manifolds in the literature, and we will use the following:

\begin{definition}A Calabi-Yau manifold is a compact K\"ahler manifold $X$ whose real first Chern class $c_1(X)\in H^2(X,\mathbb{R})$ vanishes, i.e. $c_1(X)=0$.
\end{definition}
Since $c_1(X)=-c_1(K_X)$, where $K_X$ is the canonical bundle of $X$, the Calabi-Yau condition is clearly equivalent to $K_X$ being topologically torsion. In fact we have the following result which was established independently in \cite[Theorem 3]{Bog}, \cite[Proposition 6.6]{Fuj}, \cite[Theorem 3.13]{Lie} (see also \cite[Theorem 2]{Cal} and \cite[Theorem 3]{Mat} for earlier results):
\begin{theorem}
Every Calabi-Yau manifold $X$ has $K_X$ holomorphically torsion, i.e. $K_X^{\otimes\ell}\cong\mathcal{O}_X$ for some $\ell\geq 1$.
\end{theorem}
In particular, some finite \'etale cover $\ti{X}$ of $X$ is a Calabi-Yau manifold with (holomorphically) trivial canonical bundle.

\subsection{Examples}
The classical examples of Calabi-Yau manifolds are the following:\\

\noindent (1) Complex tori, $X=\mathbb{C}^n/\Lambda$ where $\Lambda\cong\mathbb{Z}^{2n}$ is a lattice.\\

\noindent (2) Smooth projective hypersurfaces $X=\{P=0\}\subset\mathbb{CP}^{n+1}$ of degree $\deg P=n+2$.\\

\noindent (3) Generalizing the previous example, if $Y$ is a Fano $(n+1)$-fold, i.e. $K_Y^{-1}$ is ample, and $X\subset Y$ is a smooth hypersurface in the anticanonical linear system $|K_Y^{-1}|$, then $X$ is a Calabi-Yau $n$-fold.\\

\noindent (4) As another generalization of (3) we have smooth complete intersections in products of projective spaces $X^n=\{P_1=\cdots=P_m=0\}\subset \mathbb{CP}^{n_1}\times\cdots\times \mathbb{CP}^{n_k}$ where $n=\sum_{p=1}^k n_p-m$ and each $P_j$ is a multi-homogeneous polynomial of multidegree $(d_1^{(j)},\dots,d_k^{(j)})$ and we have $\sum_{j=1}^m d_p^{(j)}=n_p+1$ for all $1\leq p\leq k$.\\

\noindent (5) A simply connected Calabi-Yau surface (i.e. $\dim_{\mathbb{C}}X=2$) is called a $K3$ surface. The $2$-dimensional examples in (2) and (3) above are projective $K3$ surfaces.\\

\noindent (6) A special class of Calabi-Yau manifolds are hyperk\"ahler (or irreducible holomorphic symplectic) manifolds. These are simply connected compact K\"ahler manifolds $X$ with $\dim_{\mathbb{C}}X=2n$, with a holomorphic $2$-form $\Omega$ which is unique up to scaling and which is symplectic (i.e. it induces an isomorphism $T^{(1,0)}X\cong (T^{(1,0)}X)^*$). These are indeed Calabi-Yau since $\Omega^n$ trivializes $K_X$. Examples of hyperk\"ahler $2n$-manifolds are the Hilbert schemes $X=\mathrm{Hilb}^{[n]}(Y)$ of $n$ points on a $K3$ surface $Y$, see \cite{Be}.\\

It is a classical fact that Calabi-Yau manifolds of complex dimension $1$ are elliptic curves ($1$-dimensional tori). The Kodaira-Enriques classification of complex surfaces implies that every $2$-dimensional Calabi-Yau manifold is either a torus, a $K3$ surface, or a finite unramified quotient of these (which are called bielliptic and Enriques surfaces respectively). It is also known that there are only finitely many possible diffeomorphism types (and even complex deformation equivalence classes) for Calabi-Yau manifolds of dimension $\leq 2$, but whether this is true also in any given dimension $n\geq 3$ is an outstanding well-known open problem:

\begin{question}For any given $n\geq 3$, is there only a finite number of possible diffeomorphism types of Calabi-Yau $n$-folds?
\end{question}
One can also ask the stronger question of finitely many complex deformation equivalence classes, or the weaker questions of finitely many topological types, or even just finitely many values of the topological Euler characteristic $\chi(X)$.

\subsection{Ricci-flat K\"ahler metrics}
The fundamental result about Calabi-Yau manifolds, which gives them their name, is Yau's solution \cite{Ya} of the Calabi Conjecture \cite{Cal}:

\begin{theorem}\label{cy}
Let $X^n$ be a Calabi-Yau manifold and $\omega$ a K\"ahler metric on $X$. Then there is a unique K\"ahler metric $\ti{\omega}$ on $X$ with $[\ti{\omega}]=[\omega]$ in $H^2(X,\mathbb{R})$ and with vanishing Ricci curvature,
$$\Ric(\ti{\omega})=0.$$
\end{theorem}

More precisely, consider the Ricci curvature form $\Ric(\omega)$ of an arbitrary K\"ahler metric on $X$. This is a closed real $(1,1)$-form on $X$, locally given by $\Ric(\omega)=-\ddbar\log\det(g),$ whose cohomology class in $H^2(X,\mathbb{R})$ equals $2\pi c_1(X)$, hence it vanishes since $X$ is Calabi-Yau. By the $\de\db$-Lemma, there is a smooth function $F$ on $X$, unique up to an additive constant, such that $\Ric(\omega)=\ddbar F$. We fix this additive constant so that $$\int_X(e^F-1)\omega^n=0.$$ By the $\de\db$-Lemma again the Ricci-flat metric $\ti{\omega}$ that we seek to construct can be written as $\ti{\omega}=\omega+\ddbar\vp$ for some $\vp\in C^\infty(X,\mathbb{R})$ (unique up to an additive constant). It is then elementary to see that $\ti{\omega}$ is Ricci-flat if and only if it satisfies
$$\ti{\omega}^n=(\omega+\ddbar\vp)^n=e^F\omega^n,$$
which can be thought of as a PDE for the unknown function $\vp$, which of course should satisfy that $\omega+\ddbar\vp>0$. This is thus a complex Monge-Amp\`ere equation, a fully nonlinear (if $n\geq 2$) elliptic second order scalar PDE for $\vp$, which we can write as
\begin{equation}\label{ma1}
(\omega+\ddbar\vp)^n=e^F\omega^n,\quad \omega+\ddbar\vp>0,\quad \sup_X\vp=0,
\end{equation}
by fixing a normalization for $\vp$. Yau's fundamental result \cite{Ya} is then:

\begin{theorem}
Let $(X^n,\omega)$ be a compact K\"ahler manifold and $F\in C^\infty(X,\mathbb{R})$ be a given function with $\int_X(e^F-1)\omega^n=0.$ Then there exists a unique $\vp\in C^\infty(X,\mathbb{R})$ solving \eqref{ma1}.
\end{theorem}

The Ricci-flat K\"ahler metrics provided by Theorem \ref{cy} have restricted holonomy contained in $SU(n)$. If $X$ is a torus, then every Ricci-flat K\"ahler metric is flat and the holonomy is trivial. If $X$ is simply connected and $H^i(X,\mathcal{O}_X)=0$ for $0<i<n$ then the holonomy is exactly $SU(n)$, while if $X$ is hyperk\"ahler (of dimension $2n$) then the holonomy is exactly $Sp(n)$. The Bogomolov-Beauville-Calabi Decomposition Theorem \cite{Be,Bog2,Cal} shows that if $(X,\omega)$ is a Ricci-flat Calabi-Yau manifold, then some finite cover splits holomorphically and isometrically into the product of factors each of which is either a flat torus, a Calabi-Yau manifold $Y$ with holonomy equal to $SU(\dim Y)$ or a hyperk\"ahler manifold $Z$ with holonomy equal to $Sp(\dim Z/2)$.

\subsection{The K\"ahler cone and its boundary}
On a compact K\"ahler manifold $(X,\omega)$ the cohomology class $[\omega]$ lies naturally in the cohomology group
$H^{1,1}(X,\mathbb{R})\subset H^2(X,\mathbb{R})$ of real de Rham $2$-classes which admit a representative which is a closed real $(1,1)$-form. The vector space $H^{1,1}(X,\mathbb{R})$ has real dimension equal to the Hodge number $h^{1,1}(X)$. The cohomology classes of K\"ahler metrics form a cone
$$\mathcal{C}_X=\{[\omega]\ |\ \omega \textrm{ K\"ahler metric on }X\}\subset H^{1,1}(X,\mathbb{R}),$$
the {\em K\"ahler cone} of $X$, which is easily seen to be an open and convex cone. Its closure $\ov{\mathcal{C}_X}$ is called the {\em nef cone}, and is a closed convex and salient cone. It is easily verified that a class $[\alpha]\in H^{1,1}(X,\mathbb{R})$ is nef if and only if for every $\ve>0$ it contains a representative $\alpha_\ve$ which is a smooth real $(1,1)$-form which satisfies
\begin{equation}\label{mona}
\alpha_\ve\geq -\ve\omega,
\end{equation} on $X$.
The boundary $\de\mathcal{C}_X$ of the K\"ahler cone has in general a very complicated geometry, already on $K3$ surfaces, see \cite{Tot} for some explicit examples.

On the other hand, a class $[\alpha]$ which admits a smooth semipositive representative $\alpha\geq 0$ (i.e. satisfying \eqref{mona} with $\ve=0$) is called semipositive. The cone of semipositive classes clearly contains the K\"ahler cone and is contained in the nef cone, but this last inclusion is in general strict. The first example of a nef $(1,1)$-class with no smooth semipositive representative was constructed by Demailly-Peternell-Schneider \cite{DPS} on a ruled surface, and we will give examples on $K3$ surfaces in section \ref{below}.

\subsection{Degenerations of Ricci-flat K\"ahler metrics}
Thanks to Yau's Theorem \ref{cy}, there is a bijective correspondence between the set of Ricci-flat K\"ahler metrics on a Calabi-Yau manifold $X$ and the points of its K\"ahler cone $\mathcal{C}_X$. Because of this result, there are two basic ways in which Ricci-flat Calabi-Yau metrics can degenerate: we can
either fix the complex structure and move the K\"ahler class to the boundary of the K\"ahler cone, or we can ``fix the K\"ahler class'' and degenerate the complex structure in a family. Of course one can also consider more general degenerations in which the complex and K\"ahler structures vary simultanously, but we will not consider these here. Instead, let us explain these two basic setups more precisely.

\subsection{Setup I: K\"ahler class degenerations}\label{sect1}
We assume we have a Calabi-Yau manifold
$X$ with a cohomology class $[\alpha]\in\de\mathcal{C}_X$
and a path $[\alpha_t]\in H^{1,1}(X,\mathbb{R}), 0\leq t\leq 1,$ with $[\alpha_t]\in\mathcal{C}_X$ for $t>0$ and $[\alpha_0]=[\alpha]$.
Theorem \ref{cy} gives us a unique Ricci-flat K\"ahler metric $\omega_t$ in the class $[\alpha_t]$ for $t>0$.

In many cases, we will consider the special path $[\alpha_t]=[\alpha_0]+t[\omega_X]$, where $\omega_X$ is a K\"ahler metric on $X$.

\subsection{Setup II: complex structure degenerations}\label{sect2}
For simplicity here we restrict to the projective case. The setup is then that we have a family
\[\xymatrix{
\mathfrak{X}\ar@{^{(}->}[r] \ar[d]_\pi
&\mathbb{CP}^N\times \Delta \ar[ld]\\
\Delta &}
\]
where $\mathfrak{X}$ is an irreducible $n+1$-dimensional variety, $\Delta\subset\mathbb{C}$ is the unit disc, $\pi$ is a surjective holomorphic map with connected fibers which is a proper holomorphic submersion over $\Delta^*$ and $X_t=\pi^{-1}(t)$ are Calabi-Yau manifolds for $t\in \Delta^*$. We denote by $\mathfrak{L}\to\mathfrak{X}$ the pullback of $\mathcal{O}_{\mathbb{CP}^N}(1)$ and by $[\alpha_t]=c_1(\mathfrak{L}|_{X_t})$ for $t\in \Delta^*$. This is a K\"ahler class on $X_t$ so by Theorem \ref{cy} it contains a unique Ricci-flat K\"ahler metric $\omega_t$, for $t\in \Delta^*$.

\subsection{The main question}
The main motivating question is the following:

\begin{question}\label{mainq} What is the behaviour of the Ricci-flat manifolds $(X,\omega_t)$ in Setup I and $(X_t,\omega_t)$ in Setup II when $t$ approaches zero?
Do they have a Gromov-Hausdorff limit, or a smooth limit away from a singular set? If so, how large are the singularities limit metric space? Is the limit space homeomorphic to a complex analytic space?
\end{question}

There is of course a vast literature on the general problem in Riemannian geometry of understanding limits of Einstein manifolds, starting from the general result of Gromov \cite{gr2} which shows that there is some sequence $t_i\to 0$ such that $\omega_{t_i}$ converges in the (pointed) Gromov-Hausdorff sense to a metric space, which under certain assumptions can be proven to have additional properties, see e.g. \cite{CC1,CC2,CC3}, but even the uniqueness of the limit space in our setting does not follow from this theory. The philosophy of the work that we will overview is that in the setting of Question \ref{mainq} it is possible to go substantially beyond what the general Riemannian convergence theory gives, using tools from K\"ahler geometry and the analysis of complex Monge-Amp\`ere equations, and to give definite answers to Question \ref{mainq} in many cases.

\section{Degenerating the K\"ahler class}
In this section we will discuss Setup I from \S \ref{sect1}. Recall that we have a Calabi-Yau manifold $X^n$ with a path $[\alpha_t]\in H^{1,1}(X,\mathbb{R}), 0\leq t\leq 1,$ of classes which are K\"ahler for $t>0$ and with $[\alpha_0]=[\alpha]\in\de\mathcal{C}_X$, and for $t>0$ we let $\omega_t$ be the unique Ricci-flat K\"ahler metric in the class $[\alpha_t]$.

We can fix a unit-volume Ricci-flat K\"ahler metric $\omega_X$ on $X$, and then the Ricci-flatness of $\omega_t$ is equivalent to
\begin{equation}\label{ma2}
\omega_t^n=\left(\int_X\omega_t^n\right)\omega_X^n.
\end{equation}

\subsection{Diameter}
The first basic result about the metrics $\omega_t$ is the following:
\begin{theorem}[T. \cite{To0}, Zhang \cite{ZT}]\label{diam}
There is a constant $C>0$ such that
$$\diam(X,\omega_t)\leq C,$$
for all $0<t\leq 1$.
\end{theorem}

The idea of the proof of this result is to first use a lemma of Demailly-Peternell-Schneider \cite{DPS2} together with \eqref{ma2} (and the assumption that $[\omega_t]$ converges to the limiting class $[\alpha]$) to show that for all $t$ there is a unit-radius $\omega_t$-geodesic ball $B_t=B_{\omega_t}(x_t,1)\subset X$ with $\omega_t$-volume a definite fraction of the total volume, i.e. $\int_{B_t}\omega_t^n\geq C^{-1}\int_X\omega_t^n$. On the other hand, an inequality of Yau \cite[Theorem I.4.1]{SY} (which only uses that $\Ric(\omega_t)\geq 0$) gives that for all $0<R<\diam(X,\omega_t)$ we have
$$\frac{R-1}{4n}\leq \frac{\int_{B_{\omega_t}(x_t,2(R+1))}\omega_t^n}{\int_{B_t}\omega_t^n},$$
and taking $R=\diam(X,\omega_t)-1$ and using the result above, we obtain Theorem \ref{diam}.\\

From Theorem \ref{diam}, it follows in particular that any Gromov-Hausdorff limit of some sequence $(X,\omega_{t_i})$ ($t_i\to 0$) is a compact length metric space.

\subsection{Volume}
The volume of the Ricci-flat metrics $\omega_t$ is readily calculated as (we drop here and in the following the customary factor of $n!$)
$$\Vol(X,\omega_t)=\int_X\omega_t^n=\int_X\alpha_t^n\to\int_X\alpha^n,\quad \textrm{as }t\to 0,$$
so the behavior of the volume is completely determined by the intersection number $\int_X\alpha^n\geq 0$.

If $\int_X\alpha^n>0$ we will say that the metrics $\omega_t$ are non-collapsed, which is justified since by the Bishop-Gromov volume comparison theorem and the diameter bound in Theorem \ref{diam} we have for every $x\in X$ and $r\leq \diam(X,\omega_t)$,
$$\Vol(B_{\omega_t}(x,r),\omega_t)\geq C^{-1}r^{2n},$$
for a uniform constant $C>0$ independent of $t>0$. We will refer to the non-collapsed setting as Setup I.A.

If on the other hand we have $\int_X\alpha^n=0$, we will say that the metrics $\omega_t$ are volume collapsed. We will call this Setup I.B.

\subsection{Setup I.A: Volume non-collapsed}
In this section we discuss the volume non-collapsed case of Setup I, which we call Setup I.A: we have a Calabi-Yau manifold $X$ with a cohomology class $[\alpha]\in\de\mathcal{C}_X$ with $\int_X\alpha^n>0$ (this is known as ``nef and big''), with and a path $[\alpha_t]\in H^{1,1}(X,\mathbb{R}), 0\leq t\leq 1,$ with $[\alpha_t]\in\mathcal{C}_X$ for $t>0$ and $[\alpha_0]=[\alpha]$, and we let $\omega_t$ be the unique Ricci-flat K\"ahler metric in the class $[\alpha_t]$ for $t\neq 0$.

\begin{example}\label{kum}
The simplest example of Setup I.A is when $X$ is a Kummer $K3$ surface, associated with a $2$-torus $T=\mathbb{C}^2/\Lambda$. To construct $X$ we consider the involution $\iota:T\to T$ induced by reflection across the origin in $\mathbb{C}^2$, and take the quotient $Y=T/\iota$. Since $\iota$ has $16$ fixed points (the $2$-torsion points of $T$), $Y$ is has $16$ orbifold points (which are rational double points), and blowing up these $16$ points we get $\pi:X\to Y$ where $X$ is our Kummer $K3$ surface. We then take $\omega_Y$ any orbifold K\"ahler metric on $Y$, and let $[\alpha]=\pi^*[\omega_Y]$. This is a nef and big class on $X$ which is not K\"ahler since it intersects trivially with the $16$ exceptional divisors of $\pi$.
\end{example}

Going back to the general Setup I.A, if $V^k\subset X$ is an irreducible $k$-dimensional closed analytic subvariety, smooth real $(k,k)$-forms on $X$ can be integrated over $V$ (by restricting to the regular part) to get a finite number, which vanishes if the form is exact (see e.g. \cite{GH} for these classical results by Lelong), and in fact the $2k$-dimensional volume of $V$ with respect to $\omega_t$ can be calculated as
$$\Vol(V,\omega_t)=\int_V\omega_t^k=\int_V\alpha_t^k\to \int_V\alpha^k,\quad \textrm{as }t\to 0,$$
so again the behavior of the volume of $V$ is controlled by the intersection number $\int_V\alpha^k\geq 0$. We can collect together all the subvarieties whose volume is shrinking to zero and define the {\em null locus} of $[\alpha]$ as
$$\Null([\alpha])=\bigcup_{\int_V \alpha^{\dim V}=0}V,$$
where the (set-theoretic) union is over all positive-dimensional irreducible closed analytic subvarieties $V\subset X$ with vanishing intersection with $[\alpha]$. This object first appeared explicitly in the work of Nakamaye \cite{Nak} in the context of algebraic geometry where $X$ is a projective manifold and $[\alpha]=c_1(L)$ for some holomorphic line bundle $L\to X$.

In the Kummer Example \ref{kum} we have that $\Null([\alpha])$ equals the union of the $16$ exceptional divisors of $\pi$.

Despite it being defined as a potentially large union of subvarieties, the null locus is in fact Zariski closed as shown by Collins and the author \cite{CT}:

\begin{theorem}\label{eins}
Let $X^n$ be a compact K\"ahler manifold and $[\alpha]$ a nef $(1,1)$-class with $\int_X\alpha^n>0$. Then $\Null([\alpha])\subset X$ is a proper closed analytic subvariety, which is empty if and only if $[\alpha]\in\mathcal{C}_X$.
\end{theorem}

The fact that $\Null([\alpha])=\emptyset$ characterizes K\"ahler classes among nef classes is just a restatement of the famous Nakai-Moishezon criterion for K\"ahler classes by Demailly-P\u{a}un \cite{DP}, and the arguments in \cite{CT} together with \cite{Ch} provide a new proof of it.

Furthermore, inspired by foundational results of Nakamaye \cite{Nak} and Ein-Lazarsfeld-Musta\c{t}\u{a}-Nakamaye-Popa \cite{ELMNP} in algebraic geometry, we showed in \cite{CT} that the class $[\alpha]$ behaves like a K\"ahler class on the complement of its null locus, in the following sense:

\begin{theorem}\label{zwei}
Let $(X^n,\omega)$ be a compact K\"ahler manifold and $\alpha$ a closed real $(1,1)$ form such that $[\alpha]$ is nef with $\int_X\alpha^n>0$. Then there is a quasi-psh function $\vp$ on $X$, which is smooth on $X\backslash \Null([\alpha])$ and has logarithmic singularities along $\Null([\alpha])$, such that
$$\alpha+\ddbar\vp\geq \ve\omega,$$
in the weak sense of currents on $X$, for some $\ve>0$.
\end{theorem}
Without entering into the technicalities of closed positive currents, all that is important for us here is that outside of $\Null([\alpha])$ we have that $\alpha+\ddbar\vp$ is a smooth K\"ahler metric with a definite lower bound $\ve\omega>0$ which does not degenerate near $\Null([\alpha])$, and that the potential $\vp$ goes to $-\infty$ uniformly along $\Null([\alpha])$. One can thus think of $\alpha+\ddbar\vp$ as a K\"ahler metric with certain singularities along $\Null([\alpha])$, or in technical terms as a K\"ahler current in the class $[\alpha]$.

Using Theorems \ref{eins} and \ref{zwei}, the following result was proved in \cite{To0,RZ,CT}:

\begin{theorem}\label{noncoll}
Suppose we are in Setup I.A, so $X^n$ is Calabi-Yau and $[\alpha]\in\de\mathcal{C}_X$ is a nef $(1,1)$-class with $\int_X\alpha^n>0$. Then there is an incomplete Ricci-flat K\"ahler metric $\omega_0$ on $X\backslash\Null([\alpha])$ such that $\omega_t\to \omega_0$ in $C^\infty_{\rm loc}(X\backslash\Null([\alpha])).$ The limit $\omega_0$ is unique, independent of the path $[\alpha_t]$, and $(X,\omega_t)$ converges in the Gromov-Hausdorff topology to $(Z,d)$, the metric completion of $(X\backslash\Null([\alpha]),\omega_0)$.
\end{theorem}

This theorem applies in particular in the Kummer $K3$ surfaces of Example \ref{kum}, thus recovering earlier results of Kobayashi-Todorov \cite{KoT} and LeBrun-Singer \cite{LS}. Briefly, the idea is that $\vp$ can be used as a barrier function to obtain uniform $C^\infty$ estimates for $\omega_t$ on compact subsets away from $\Null([\alpha])$. In the case when $[\alpha]\in H^2(X,\mathbb{Q})$ one can use Nakamaye's Theorem \cite{Nak} in place of Theorems \ref{eins} and \ref{zwei}, which was done by the author in \cite{To0}, and the Gromov-Hausdorff limit in this case was obtained by Rong-Zhang \cite{RZ}. Lastly, Collins and the author \cite{CT} used Theorems \ref{eins} and \ref{zwei} to remove the rationality assumption.

Furthermore, if we assume that $[\alpha]\in H^2(X,\mathbb{Q})$, then some multiple of $[\alpha]$ equals $c_1(L)$ for a holomorphic line bundle $L\to X$, and the fact that $c_1(L)$ is nef and big implies that $X$ is projective (since it is K\"ahler and Moishezon, see e.g. \cite{DP}) and Kawamata's base-point-free theorem \cite{KMM} implies that $L$ is semiample, so the linear system induced by some high tensor power of $L$ defines a surjective holomorphic birational map $\pi:X\to Y$ with connected fibers onto a normal projective variety $Y$ which is a Calabi-Yau variety with at worst canonical singularities (see e.g. \cite{To0}). Furthermore, the exceptional locus of $\pi$ is equal to $\Null([\alpha])$, we have that $L=f^*A$ for some ample line bundle $A$ on $Y$, and the limiting Ricci-flat metric $\omega_0$ is naturally the pullback of a singular Ricci-flat metric on $Y$ \cite{EGZ} in a class proportional to $c_1(A).$  In this case it was proved by Song \cite{So} (improving results of Donaldson-Sun \cite{DS} that do not apply here because of the changing K\"ahler classes) that the Gromov-Hausdorff limit $Z$ is in fact homeomorphic to $Y$, and so the non-collapsed Gromov-Hausdorff limit of $(X,\omega_t)$ has an algebraic structure. Furthermore, the singular set of $Z$ (in the sense of \cite{CC1}) has real Hausdorff codimension at least $4$. 

These results give an answer to Question \ref{mainq} in Setup I.A when $[\alpha]$ is rational. We expect that a similar picture should hold even when the class $[\alpha]$ (or any of his multiples $\lambda[\alpha],\lambda>0$) is not rational:
\begin{conjecture}\label{1}
Let $X^n$ be a Calabi-Yau manifold and $[\alpha]\in\de\mathcal{C}_X$ a nef $(1,1)$-class with $\int_X\alpha^n>0$. Then there is a surjective holomorphic and bimeromorphic map $\pi:X\to Y$ with connected fibers onto a normal compact K\"ahler analytic space $Y$ such that $[\alpha]=\pi^*[\beta]$ for some K\"ahler class $[\beta]$ on $Y$ (in the sense of e.g. \cite{EGZ}), and $\mathrm{Exc}(\pi)=\Null([\alpha])$.
\end{conjecture}

In particular, the nef and big class $[\alpha]$ would contain a smooth semipositive representative, which is also an open problem. The following related conjecture is also open:
\begin{conjecture}\label{2}
Let $X^n$ be a Calabi-Yau manifold and $\alpha$ a closed real $(1,1)$-form such that $[\alpha]\in\de\mathcal{C}_X$ is nef. Then there is a bounded function $\vp\in L^\infty(X)$ such that $\alpha+\ddbar\vp\geq 0$ in the weak sense of currents on $X$.
\end{conjecture}

As we just said, if Conjecture \ref{1} holds then we can even find such a $\vp\in C^\infty(X)$ when $\int_X\alpha^n>0$, so Conjecture \ref{2} when $[\alpha]$ is nef and big follows from Conjecture \ref{1}. Conjecture \ref{1} is known when $n=2$ by work of Filip and the author \cite{FT} and when $n=3$ by H\"oring \cite{Hor}. Note also that if $[\alpha]\in H^2(X,\mathbb{Q})$ then every irreducible component of $\Null([\alpha])=\mathrm{Exc}(\pi)$ is uniruled, as follows for example from \cite{Kaw}, and it is natural to conjecture that this property also holds for general nef and big $[\alpha]$.

On the other hand, Conjecture \ref{2} when $\int_X\alpha^n=0$ is not even known for $K3$ surfaces, see also section \ref{later}. Observe that the condition in Conjecture \ref{2} is equivalent to the fact that a closed positive current with minimal singularities in the class $[\alpha]$ has bounded potentials, in the language of Demailly \cite{DPS}. The Calabi-Yau condition is crucial in both conjectures, see e.g. \cite{DPS} for a counterexample to Conjecture \ref{2} where $X$ is a ruled surface. Also, as remarked above, it would follow from Conjecture \ref{1} that every nef and big class $[\alpha]$ on a Calabi-Yau manifold admits a smooth semipositive representative, but this is known to be false for nef classes with $\int_X\alpha^n=0$ as discussed in detail in section \ref{below} below.

Going back to the Ricci-flat metrics $\omega_t$, in analogy with Theorem \ref{noncoll} and the aforementioned result of \cite{So} we then expect the following:
\begin{conjecture}\label{3}
In setup I.A, the Gromov-Hausdorff limit $Z$ of the Ricci-flat manifolds $(X,\omega_t)$ is homeomorphic to a normal compact K\"ahler analytic space $Y$. The singular set of $Z$ (in the sense of \cite{CC1}) has real Hausdorff codimension at least $4$. If Conjecture \ref{1} holds, $Y$ should be the same space that appears there.
\end{conjecture}

Lastly, let us mention that one can ask further questions about the structure of the Gromov-Hausdorff limit space $(Z,d)$ in Theorem \ref{noncoll} near its singular points. In particular, it is not known whether the tangent cones at singular points of $(Z,d)$ are unique. This is established in \cite{DS2} in a related but different situation (with fixed polarization), and in \cite{HS} in the case when the contracted space $Y$ has certain types of isolated singularities.

\subsection{Setup I.B: volume collapsed}
In this section we discuss the volume collapsed case of Setup I, which we call Setup I.B: we have a Calabi-Yau manifold $X^n$ with a cohomology class $[\alpha]\in\de\mathcal{C}_X$ with $\int_X\alpha^n=0$, which we approach with the path $[\alpha_t]=[\alpha]+t[\omega_X]$ (where $\omega_X$ is a K\"ahler metric on $X$), and we let $\omega_t$ be the unique Ricci-flat K\"ahler metric in the class $[\alpha_t]$ for $t>0$.

The choice of path $[\alpha]+t[\omega_X]$ comes from the work of Gross-Wilson \cite{GW}, and it is also formally analogous to the path traced in cohomology by the K\"ahler-Ricci flow on compact K\"ahler manifolds with semiample canonical bundle and intermediate Kodaira dimension \cite{ST} (in fact, many of the results described below have analogs in the setting of the K\"ahler-Ricci flow, see \cite{To4}). 

\begin{example}\label{k3}
The simplest example of Setup I.B is when $X$ is a $K3$ surface that admits an elliptic fibration $f:X\to \mathbb{CP}^1$. This is a surjective holomorphic map with connected fibers and with all but finitely many fibers equal to elliptic curves. If we take $\omega_{\rm FS}$ to be the Fubini-Study metric on $\mathbb{CP}^1$, then $[\alpha]=f^*[\omega_{\rm FS}]$ is a nef class on $X$ with $\int_X\alpha^2=0$. The behavior of the Ricci-flat metrics $\omega_t$ on $X$ in this case was first investigated by Gross-Wilson \cite{GW}.
\end{example}

\subsection{Fiber spaces}
Example \ref{k3} can be generalized to higher dimension as follows. Let $X^n$ be a Calabi-Yau manifold and suppose that we have $f:X\to Y$ a surjective holomorphic map with connected fibers onto a compact K\"ahler analytic space $Y$ with $0<\dim Y=m<\dim X$. The map $f$ will be called a {\em fiber space}. If $[\omega_Y]$ is a K\"ahler metric on $Y$ (in the sense of analytic spaces \cite{EGZ}) then $[\alpha]=f^*[\omega_Y]$ is a nef class on $X$ with $\int_X\alpha^n=0$.

Many nef classes $[\alpha]$ arise in this way: for example if $X$ is a $K3$ surface and $0\neq [\alpha]\in \de\mathcal{C}_X$ satisfies $\int_X\alpha^2=0$ and $[\alpha]\in H^2(X,\mathbb{Q})$ then $[\alpha]$ is the pullback of a K\"ahler class from the base $\mathbb{CP}^1$ of an elliptic fibration on $X$, see \cite{FT}. In higher dimensions, a well-known conjecture in algebraic geometry (see e.g. \cite{LOP}) would imply that if $X$ is any projective Calabi-Yau manifold and $0\neq [\alpha]\in \de\mathcal{C}_X$ satisfies $\int_X\alpha^n=0$ and $[\alpha]\in H^2(X,\mathbb{Q})$ then $X$ admits a fiber space structure $f:X\to Y$ and $[\alpha]=f^*[\omega_Y]$ for some K\"ahler class $[\omega_Y]$ on $Y$. On the other hand, the analogous conjecture for $[\alpha]=c_1(D)$ with $D$ a nef $\mathbb{R}$-divisor fails, see Theorem \ref{dyn}, so the rationality assumption is crucial here.

It is important to note that fiber spaces in general may have singular fibers. More precisely, if we denote by $S'\subset Y$ the union of the singular locus of $Y$ together with the critical values of $f$ on the regular part of $Y$, then $S'$ is a closed analytic subvariety of $Y$, its preimage $S=f^{-1}(S')$ is a closed analytic subvariety of $X$, and $f:X\backslash S\to Y\backslash S'$ is a proper holomorphic submersion with fibers $X_y=f^{-1}(y)$ Calabi-Yau manifolds of dimension $n-m$. The subvariety $S$ is informally referred to as the union of singular fibers of $f$.

It is also important to note that the smooth fibers $X_y=f^{-1}(y)$ for $y\in Y\backslash S'$ are all pairwise diffeomorphic (Ehresmann's Lemma) but in general are not biholomorphic. By a result of Fischer-Grauert \cite{FG}, the fibers $X_y$ are all pairwise biholomorphic for all $y$ in some open subset $U\subset  Y\backslash S'$ if and only if $f:f^{-1}(U)\to U$ is a holomorphic fiber bundle.

When $Y$ is smooth and $f$ is a submersion everywhere (i.e. with our notation when $S=\emptyset$) then we have the following result proved by Zhang and the author \cite{TZ0,TZ2}:
\begin{theorem}\label{triv}
Let $X$ be a compact Calabi-Yau manifold and $f:X\to Y$ a holomorphic submersion with connected fibers onto a lower-dimensional K\"ahler manifold. Then $f$ is a holomorphic fiber bundle, and both $Y$ and the fiber of $f$ are Calabi-Yau.
\end{theorem}

Furthermore, if $X$ is projective, then this fiber bundle is actually the quotient of the trivial bundle by a finite group action. This result shows that to have nontrivial variation of the complex structure of the fibers of $f$, one must allow singular fibers.

\subsection{Setup I.B.1: fiber spaces}
The discussion of our main Question \ref{mainq} in the collapsing Setup I.B now splits into two subcases, according to whether the class $[\alpha]$ comes from the base of a fiber space structure on $X$, or not. We first discuss in detail the case when it does, so we have $f:X\to Y$ a fiber space with $[\alpha]=f^*[\omega_Y]$ for some K\"ahler class $[\omega_Y]$ on $Y$, and we will call this Setup I.B.1.

The Ricci-flat K\"ahler metrics $\omega_t$ in the class $[\alpha_t]=f^*[\omega_Y]+t[\omega_X]$ can be written as
$$\omega_t=f^*\omega_Y+t\omega_X+\ddbar\vp_t,\quad \sup_X\vp_t=0,$$
where the normalized K\"ahler potentials $\vp_t$ are uniquely determined and they solve the degenerating family of complex Monge-Amp\`ere equations
\begin{equation}\label{ma4}
(f^*\omega_Y+t\omega_X+\ddbar\vp_t)^n=c_t t^{n-m}\omega_X^n,
\end{equation}
where $c_t$ is a constant which bounded away from $0$ and $\infty$.

\subsubsection{Behavior away from the singular fibers}
The first general result about \eqref{ma4} is a uniform bound for $\vp_t$,
$$\|\vp_t\|_{L^\infty(X)}\leq C,$$
which was proved independently in \cite{DPa,EGZ2}, extending classical work of Ko\l odziej \cite{Kol}. Using this estimate as a starting point, and using also Yau's Schwarz Lemma computations \cite{Ya2}, the author proved in \cite{To1} (extending \cite{ST0} which dealt with elliptically fibered $K3$ surfaces):

\begin{theorem}\label{col}
Suppose we are in Setup I.B.1, so $X^n$ is Calabi-Yau, $f:X\to Y$ is a fiber space, $[\alpha]=f^*[\omega_Y]$ as above, $[\alpha_t]=[\alpha]+t[\omega_X]$ and $\omega_t$ is the Ricci-flat K\"ahler metric on $X$ in the class $[\alpha_t]$ for $t>0$. Then there is a K\"ahler metric $\omega_0=\omega_Y+\ddbar\vp_0$ on $Y\backslash S'$ with
$$\Ric(\omega_0)=\omega_{\rm WP}\geq 0,$$
such that $\omega_t\to f^*\omega_0$ as $t\to 0$ in the sense that $\vp_t\to f^*\vp_0$ in  $C^{1,\gamma}_{\rm loc}(X\backslash S,\omega_X)$ for every $\gamma<1$. Furthermore, on every compact subset of $X\backslash S$ we have
\begin{equation}\label{lastima}
C^{-1}(f^*\omega_Y+t\omega_X)\leq \omega_t\leq C(f^*\omega_Y+t\omega_X).
\end{equation}
for a uniform constant $C>0$.
\end{theorem}
Thanks to \eqref{lastima}, the Ricci-flat metrics $\omega_t$ shrink the smooth fibers $X_y$ of $f$ to points as $t\to 0$ and collapse to the metric $\omega_0$ on the base in the $C^{1,\gamma}$ topology of K\"ahler potentials. The Weil-Petersson form $\omega_{\rm WP}$ measures the variation of the complex structure of the smooth fibers $X_y$, and in particular it vanishes whenever $f$ is a holomorphic fiber bundle, and already appears in \cite{GW}.

It is highly desirable to improve the topology of the convergence of $\omega_t$ to $f^*\omega_0$ (away from the singular fibers) and the main open question is then the following (cf. \cite{To2,To3}):
\begin{conjecture}\label{smooth}
In Setup I.B.1, we have that $\omega_t\to f^*\omega_0$ in $C^\infty_{\rm loc}(X\backslash S,\omega_X)$.
\end{conjecture}
In other words, one has to establish a uniform bound for $\|\vp_t\|_{C^k(K,\omega_X)}$ on every compact subset $K\subset X\backslash S$ for every $k\in\mathbb{N}$. It follows from the work of Fine \cite{Fi} that Conjecture \ref{smooth} holds when $S=\emptyset$, but recall however that thanks to Theorem \ref{triv} this assumption is very restrictive as it implies that $f$ is a holomorphic fiber bundle. When singular fibers are present, the currently known results are:

\begin{theorem}\label{co1}
Conjecture \ref{smooth} holds in any of the following cases:
\begin{itemize}
\item (Gross-T.-Zhang \cite{GTZ}, Hein-T. \cite{HT}) if the smooth fibers $X_y$ are tori (or finite \'etale quotients of tori). In this case the metrics $\omega_t$  have locally uniformly bounded sectional curvature away from $S$.
\item (Hein-T. \cite{HT2}) if the smooth fibers $X_y$ are pairwise biholomorphic (but $S$ may be nonempty)
\end{itemize}
\end{theorem}

We will give some ideas of the proof of this theorem in section \ref{ideas} below. On the other hand, without any extra assumptions on the fibers of $f$, the convergence in Theorem \ref{col} was improved by Weinkove, Yang and the author \cite{TWY} to $C^0_{\rm loc}(X\backslash S,\omega_X)$ convergence of $\omega_t$ to $f^*\omega_0$. This was more recently further improved by Hein and the author \cite{HT2}:
\begin{theorem}\label{co3}
In Setup I.B.1, we have that $\omega_t\to f^*\omega_0$ in $C^\gamma_{\rm loc}(X\backslash S,\omega_X)$ for every $\gamma<1$.
\end{theorem}
Theorems \ref{co1} and \ref{co3} are in fact purely local on the base.

\subsubsection{Adiabatic limits}
The above results show that the Ricci-flat metrics $\omega_t$ shrink the smooth fibers and collapse to a metric on the base. The second result that we will discuss, proved by Zhang and the author \cite{TZ}, shows that if we ``zoom in'' to make the fibers have fixed size (and therefore the base of the fibration is getting blown up), then we see a ``semi-Ricci-flat'' picture. This is known as an ``adiabatic limit''.

\begin{theorem}\label{co2}
In Setup I.B.1, given any $x\in X\backslash S$ fix a neighborhood $U$ of $y=f(x)$ relatively compact in $Y\backslash S'$, Then the pointed spaces $(f^{-1}(U),\frac{\omega_t}{t},x)$ converge smoothly modulo diffeomorphisms to the cylinder $$(\mathbb{C}^m\times X_y, \omega_{\mathbb{C}^m}\oplus\omega_{{\rm SRF},y}, p),$$ where $\omega_{\mathbb{C}^m}$ is a Euclidean metric and $\omega_{{\rm SRF},y}$ denotes the unique Ricci-flat K\"ahler metric on $X_y$ cohomologous to $[\omega_X|_{X_y}]$.
\end{theorem}

We emphasize here that the complex structure on $f^{-1}(U)$ also converges smoothly (modulo diffeomorphisms) to the product complex structure of $\mathbb{C}^m\times X_y$. Furthermore, the diffeomorphisms that appear in this theorem are completely explicit, and are obtained by fixing a smooth trivialization of $f$ (up to shrinking $U$) and stretching out the base directions while being the identity on the fibers.

As a consequence of Theorem \ref{co2}, if the smooth fibers $X_y$ are not tori or quotients, then the sectional curvatures of $\omega_t$ have to blow up (to $+\infty$ and also to $-\infty$) on every fiber $X_y$, since the rescaled metrics $\frac{\omega_t}{t}$ converge smoothly modulo diffeomorphism to the Ricci-flat (but non-flat) cylinder as above which has both positive and negative sectional curvatures in the fiber directions.

Of course, all of these results were already obtained by Gross-Wilson \cite{GW} using their gluing construction for elliptically fibered $K3$ surfaces with only $I_1$ singular fibers. More recently, Li \cite{Li,Li2} proves a similar result also when we zoom in at $x$ on a singular fiber when $X$ is a Calabi-Yau $3$-fold, $X_y$ are $K3$ surfaces, and $f$ is a Lefschetz fibration.

\subsubsection{Gromov-Hausdorff limits}
Finally, we discuss the Gromov-Hausdorff limit of $(X,\omega_t)$. Recall that thanks to the diameter bound in Theorem \ref{diam} and to Gromov's precompactness theorem, for every sequence $t_i\to 0$ we can find a subsequence (still denoted by $t_i$) and a compact metric space $(Z,d)$ such that $(X,\omega_{t_i})\to (Z,d)$ in the Gromov-Hausdorff topology, as $i\to\infty$.

It is then natural to ask whether $(Z,d)$ depends on the chosen sequence or is in fact unique, and what its relation is to the collapsed limit $(Y\backslash S',\omega_0)$ away from the singular fibers given by Theorem \ref{co1}.

The first general result in this direction follows by combining the work of Gross, Zhang and the author \cite{GTZ} and of Weinkove, Yang and the author \cite{TWY}:

\begin{theorem}
In Setup I.B.1, suppose that $(X,\omega_{t_i})\to (Z,d)$ in the Gromov-Hausdorff topology. Then there is a map $F:(Y\backslash S',\omega_0)\to (Z,d)$ which is a homeomorphism with its image $Z_0\subset Z$, an open dense subset of $Z$, and $F$ is a local isometry.
\end{theorem}
This shows that any collapsed Gromov-Hausdorff limit $(Z,d)$ is in fact a smooth manifold on a dense open set.

The main outstanding question in then the following \cite{To2,To3}:

\begin{conjecture}\label{gh}
In Setup I.B.1, let $(Z,d)$ be the metric completion of $(Y\backslash S',\omega_0)$, and let $S_Z=Z\backslash (Y\backslash S')\subset Z$. Then
\begin{itemize}
\item[(a)] $(Z,d)$ is a compact metric space and $(X,\omega_t)\to (Z,d)$ in the Gromov-Hausdorff topology as $t\to 0$
\item[(b)] $S_Z$ has real Hausdorff codimension at least $2$
\item[(c)] $Z$ is homeomorphic to $Y$
\end{itemize}
\end{conjecture}

The conjecture is known in several cases:

\begin{theorem}
Conjecture \ref{gh} holds in any of the following cases:
\begin{itemize}
\item (Gross-T.-Zhang \cite{GTZ,GTZ2}, T.-Weinkove-Yang \cite{TWY}) if $X$ is projective and $\dim Y=1$
\item (T.-Zhang \cite{TZ2} if $X$ is projective and hyperk\"ahler
\item (Gross-T.-Zhang \cite{GTZ3}) if $X$ is projective, $Y$ is smooth and the divisorial components of $S'$ have simple normal crossings
\end{itemize}
Furthermore,
\begin{itemize}
\item (Song-Tian-Zhang \cite{STZ}) if $X$ is projective then part (a) holds in general; if furthermore $Y$ is smooth then part (c) holds
\end{itemize}
\end{theorem}

While, as discovered in \cite{STZ}, the proof of part (a) does not require a detailed understanding of the blowup behavior of the limiting metric $\omega_0$ near $S'$, so far this seems to be needed to prove part (b), and this is what is achieved in \cite{GTZ2,TZ2,GTZ3} under the above assumptions. More specifically, under the assumption that $X$ is projective, $Y$ is smooth and the divisorial components of $S'$ have simple normal crossings, it is proved in \cite{GTZ3} that near these components $\omega_0$ is quasi-isometric to a {\em conical metric} (with rational cone angles in $(0,2\pi]$) up to logarithmically blowing up errors (which are negligible compared to the conical singularities), and this is crucial to prove part (b) of Conjecture \ref{gh}.

It is also interesting to note that in the hyperk\"ahler case, in \cite{TZ2} it is proved that the Gromov-Hausdorff limit $Z$ is homeomorphic to the base $Y$ (without assuming that $Y$ is smooth), while a conjecture in algebraic geometry predicts that the base of every nontrivial fiber space $f:X\to Y$ with $X^{2n}$ projective hyperk\"ahler should be isomorphic to $\mathbb{CP}^{n}$, but this is only known in general when $Y$ is smooth \cite{Hw}.

\subsection{Smooth collapsing}\label{ideas}
In this section we give some ideas of the proof of Theorem \ref{co1} that gives smooth collapsing (away from the singular fibers) in Setup I.B.1 when the fibers are either tori or pairwise biholomorphic.

\subsubsection{Torus fibers}
First, we discuss the easier case of torus fibers \cite{GTZ,HT}. In this case, by old results in deformation theory \cite{We}, the preimage $f^{-1}(U)$ of every sufficiently small ball $U\subset Y\backslash S'$ has universal cover biholomorphic to $U\times\mathbb{C}^{n-m}$, and we denote $p:U\times\mathbb{C}^{n-m}\to f^{-1}(U)$, which satisfies $f\circ p=\mathrm{pr}_U$. In \cite{GTZ,He,HT} one then constructs a holomorphic period map for the torus fibers (equipped with a K\"ahler class from the total space) $Z:U\to \mathfrak{H}_{n-m}$ into the Siegel upper half space of symmetric $(n-m)\times (n-m)$ complex matrices with positive definite imaginary part. This is used to define an explicit semi-flat form $\omega_{\rm SF}\geq 0$ on $f^{-1}(U)$ by defining its pullback via $p$ to be
$$p^*\omega_{\rm SF}=-\frac{\mn}{4}\de\db\left(\sum_{i,j=1}^{n-m}(\mathrm{Im}Z(y))^{-1}_{ij}(z_i-\ov{z}_i)(z_j-\ov{z}_j)\right),$$
and checking that it descends to $f^{-1}(U)$, and satisfies the crucial scaling property that
$$\lambda_t^*p^*\omega_{\rm SF}=t^{-1}p^*\omega_{\rm SF},\quad \lambda_t(y,z)=(y,t^{-\frac{1}{2}}z).$$
The semi-flat form $\omega_{\rm SF}$ (in the projective case) was used in the context of elliptically fibered $K3$ surfaces by Greene-Shapere-Vafa-Yau \cite{GSVY}, and can be traced back to Satake's book \cite[Lemma IV.8.5]{Sat}.

Thanks to \eqref{lastima} it follows easily that on $f^{-1}(U)$ we have
$$C^{-1}(f^*\omega_Y+t\omega_{\rm SF})\leq \omega_t\leq C(f^*\omega_Y+t\omega_{\rm SF}),$$
and using the scaling property of $\omega_{\rm SF}$ this implies that
$$C^{-1}p^*(f^*\omega_Y+\omega_{\rm SF})\leq \lambda_t^*p^*\omega_t\leq Cp^*(f^*\omega_Y+\omega_{\rm SF}),$$
holds on $U\times\mathbb{C}^{n-m}$, and so $\lambda_t^*p^*\omega_t$ is uniformly equivalent to the Euclidean metric on any fixed compact set. Standard higher order estimates (see e.g. \cite{HT2}) then give locally uniform $C^k$ bounds on $\lambda_t^*p^*\omega_t$, for all $k$, from which it is elementary to deduce uniform $C^k$ bounds for $\omega_t$ on $f^{-1}(U)$, as well as a uniform sectional curvature bound, as desired. The appendix of \cite{DJZ} also provides an explicit rate of convergence.

\subsubsection{Isomorphic fibers} Secondly, we discuss the second case in Theorem \ref{co1}, where the smooth fibers of $f$ are fiberwise biholomorphic, following \cite{HT2}. By the Fischer-Grauert theorem, $f$ is a holomorphic fiber bundle away from the singular fibers $S$, so $f$ can be holomorphically trivialized over any sufficiently small ball $U\subset Y\backslash S'$, so in this trivialization the fiber space equals the projection $U\times F\to U$ where $F=X_y$ is a fixed Calabi-Yau $(n-m)$-fold. We let $\omega_F$ be the Calabi-Yau metric on $F$ in the class that corresponds to $[\omega_X|_{X_y}]$ under our trivialization,  let $\omega_{\mathbb{C}^m}$ be a Euclidean metric on $U$, and define a family of shrinking product K\"ahler metrics $\ti{\omega}_t=\omega_{\mathbb{C}^m}+t\omega_F$ on $U\times F$.

After a gauge-fixing argument (see \cite{HT2}) we may assume that our Ricci-flat metrics $\omega_t$ can be written as $\omega_t=\ti{\omega}_t+\ddbar\vp_t$, they satisfy the complex Monge-Amp\`ere equation
\begin{equation}\label{ma5}
\omega_t^n=(\ti{\omega}_t+\ddbar\vp_t)^n=c_t t^{n-m}\omega_{\mathbb{C}^m}^m\wedge\omega_F^{n-m},
\end{equation}
and from \eqref{lastima}
\begin{equation}\label{lastimaz}
C^{-1}\ti{\omega}_t\leq \omega_t\leq C\ti{\omega}_t,
\end{equation}
on $U\times F$.

The goal would be to prove that
\begin{equation}\label{birsa}
\|\omega_t\|_{C^k(U'\times F,\omega_{\mathbb{C}^m}+\omega_F)}\leq C,
\end{equation}
on a smaller relatively compact subdomain $U'\subset U$ ($C$ here depends on $k$, on the subdomain $U'$ and on the background data), but in fact we will prove the much stronger estimate
\begin{equation}\label{bersa}
\|\omega_t\|_{C^k(U'\times F,\ti{\omega}_t)}\leq C.
\end{equation}
Note that in the case of torus fibers, the argument above does indeed prove such a stronger estimate (on $U\times\mathbb{C}^{n-m}$) (see also \cite{JS}).

However, when the fibers are not tori, we lack the ``fiberwise scaling trick'' above (using $\lambda_t$). We also cannot resort to purely local estimates for the Monge-Amp\`ere equation \eqref{ma5} since there are simple counterexamples to local higher-order interior estimates for the analogous Monge-Amp\`ere equation on a polydisc in $\mathbb{C}^{n-m}\times\mathbb{C}^m$ satisfying \eqref{lastimaz}, see the introduction of \cite{HT2}. These estimates are also false in a global setting where the fibers of $f$ are manifolds with boundary (and Dirichlet conditions are imposed), see again \cite{HT2}. So our argument to prove \eqref{bersa} will exploit the fact that the fiber $F$ is compact without boundary.

The uniform estimate \eqref{lastimaz} means that in the collapsed geometry of $\ti{\omega}_t$ the PDE \eqref{ma5} is uniformly elliptic. Nevertheless, \eqref{bersa}  is not a consequence of any of the existing  techniques to obtain higher-order estimates for complex Monge-Amp\`ere equations, since these reference metrics are collapsing in the fiber directions. Instead, our arguments rely on a contradiction and blowup strategy, employing linear and nonlinear Liouville theorems and new sharp Schauder estimates on cylinders. The $C^\gamma$ bound in Theorem \ref{co3} is proved in a similar fashion, but with formidable extra difficulties caused by the varying complex structure of the fibers $X_y$.

Let us discuss the proof of \eqref{bersa} for $k=1$ (the estimate for $k=0$ is given by \eqref{lastimaz}). We shall identify $U$ with the Euclidean unit ball $B_1$. Ignoring issues with the boundary of $U$ and with having to replace $t$ by some sequence $t_i\to 0$, we will assume for a contradiction that
$$\mu_t=\sup_{U\times F}|\nabla^{\ti{\omega}_t}\omega_t|_{\ti{\omega}_t}\to+\infty,$$
define a stretching map in the base directions $\Psi_t:B_{\mu_t}\times F\to B_1\times F$ by $\Psi_t(z,y)=(\mu_t^{-1}z,y)$, and let replacing $\omega_t$ by $\mu_t^2\Psi_t^*\omega_t$ (still denoting this by $\omega_t$) and $\ti{\omega}_t$ by $\mu_t^*\Psi_t^*\ti{\omega}_t=\omega_{\mathbb{C}^m}+\delta_t^2\omega_F, \delta_t^2=t\mu_t^2$ (still denoting this by $\ti{\omega}_t$), we have the same estimate \eqref{lastimaz} on $B_{\mu_t}\times F$ but now $\sup_{B_{\mu_t}\times F}|\nabla^{\ti{\omega}_t}\omega_t|_{\ti{\omega}_t}=1$, achieved at $(0,y_t)$ say.

Up to passing to subsequences there are $3$ cases. In case 1 we have $\delta_t\to\infty$. Here we fix some chart on $F$ containing all the $y_t$'s for $t$ small, and in this chart times $B_{\mu_t}$ we pull back everything via the diffeomorphism $(z,y)\mapsto(z,y_t+\delta_t^{-1}y)$, so that $\omega_t$ pulls back to a Ricci-flat metric uniformly equivalent to Euclidean on $B_{\mu_t}\times B_{\delta_t}$, so standard higher order estimates give uniform $C^\infty$ bounds for $\omega_t$, which will subconverge smoothly to a Ricci-flat K\"ahler metric $\omega_0$ on $\mathbb{C}^n$, uniformly equivalent to the Euclidean metric $\omega_{\mathbb{C}^n}$ with $|\nabla^{\omega_{\mathbb{C}^n}}\omega_0|(0)=1$. Such a metric cannot exist by the Liouville Theorem for the complex Monge-Amp\`ere equation in $\mathbb{C}^n$ (see \cite{HT2}).

In case 2 (up to scaling) we have $\delta_t\to 1$. Here we argue analogously to case 1, without applying any further diffeomorphism since $\omega_t$ is already uniformly equivalent to Euclidean locally, so $\omega_t$ will subconverge smoothly to a Ricci-flat K\"ahler metric on $\mathbb{C}^m\times F$ which is uniformly equivalent (and cohomologous) to $\omega_{\mathbb{C}^m}+\omega_F$ and satisfies $|\nabla^{\omega_{\mathbb{C}^m}+\omega_F}\omega_0|(0,y_0)=1$. Such a metric cannot exist by the Liouville Theorem for the complex Monge-Amp\`ere equation on the cylinder $\mathbb{C}^m\times F$ of Hein \cite{He2} (see also \cite{LLZ}).

We are left with case 3 where $\delta_t\to 0$, which is where all the work goes. By construction we have $$\|\omega_t\|_{C^1(B_{\mu_t}\times F,\ti{\omega}_t)}\leq C.$$ We first need to improve this to \begin{equation}\label{bega}
\|\omega_t\|_{C^{1,\gamma}(K,\ti{\omega}_t)}\leq C,
\end{equation} for any fixed $\gamma<1$ and compact set $B\subset B_{\mu_t}\times F$, which is done by applying the following new Schauder estimates on an infinite cylinder, proved in \cite{HT2}:
\begin{theorem}\label{scha}
Let $(F,g_F)$ be a closed Riemannian manifold and $g_P=g_{\mathbb{R}^d}+g_F$ the product metric on the cylinder $\mathbb{R}^d\times F$. For every $k\geq 2, 0<\gamma<1,$ there is $C>0$ such that for all $0<\rho<R$, all basepoints in $\mathbb{R}^d\times F$ and all $f\in C^{k,\gamma}_{\rm loc}(B_{2R})$ we have
$$[\nabla^k f]_{C^\gamma(B_\rho)}\leq C\left([\nabla^{k-2}\Delta f]_{C^\gamma(B_R)}+(R-\rho)^{-k-\gamma}\|f\|_{L^\infty(B_R)}\right),$$
where here geodesic balls and differential operators are those of $g_P$.
\end{theorem} 
This estimate is proved using the method of Simon \cite{LSi}, by reducing to a linear Liouville Theorem (for harmonic functions) after contradiction and blowup, but unlike in \cite{LSi} (which worked in $\mathbb{R}^d$) here we also encounter $3$ cases which blow up to $\mathbb{R}^{d+\dim F}, \mathbb{R}^d\times F,$ and $\mathbb{R}^d$ respectively.

To prove \eqref{bega} the idea is then to apply the further stretching $(z,y)\mapsto (\delta_t z,y)$ and scaling the resulting metrics by $\delta_t^{-2}$, so that $\ti{\omega}_t$ becomes a fixed product metric $\omega_P$ and $\omega_t$ is uniformly equivalent to it, so by standard higher order estimates it will subconverge smoothly to some K\"ahler metric on $\mathbb{C}^m\times F$, which can be assumed with loss of generality to be equal to $\omega_P$. In this stretched picture we can then effectively linearize the Monge-Amp\`ere equation at $\omega_P$ and apply Theorem \ref{scha} to obtain the right estimate which after un-stretching and un-scaling gives \eqref{bega}.
 
Now \eqref{bega} implies in particular a $C^{1,\gamma}$ bound with respect to any fixed metric, so Ascoli-Arzel\`a gives that, up to passing to a sequence, $\omega_t\to\omega_0$ in $C^{1,\gamma}_{\rm loc}$ where $\omega_0$ is the pullback of a $C^{1,\gamma}$ K\"ahler metric on $\mathbb{C}^m$, uniformly equivalent to Euclidean.

We then argue that $\sup_{\mathbb{C}^m}|\nabla^{\omega_{\mathbb{C}^m}}\omega_0|=1$, which however is not a trivial consequence of $|\nabla^{\ti{\omega}_t}\omega_t|_{\ti{\omega}_t}(0,y_t)=1$, since the reference metrics $\ti{\omega}_t$ are collapsing in the fiber directions. Schematically, in local product coordinates, we have at $(0,y_t)$
$$1=|\nabla^{\ti{\omega}_t}\omega_t|_{\ti{\omega}_t}\sim |\nabla_z g_{zz}|+\delta_t^{-1}|\nabla_z g_{zy}|+\cdots+\delta_t^{-3}|\nabla_y g_{yy}|,$$
where we write $g_{\bullet\bullet}$ for the components of $\omega_t$,
and the claim is then that all the terms on the RHS except the first one must go to zero, which then passing to the limit gives us $\sup_{\mathbb{C}^m}|\nabla^{\omega_{\mathbb{C}^m}}\omega_0|=1$. To see why these go to zero we use the following elementary result (which uses crucially that $F$ is compact without boundary): given a metric vector bundle $E\to F$, and fix a metric on $F$, then for every $\gamma<1$ there is $C>0$ such that for every smooth section $\sigma$ of $E$ over $F$ we have $$\|\nabla\sigma\|_{L^\infty(F)}\leq C[\nabla\sigma]_{C^\gamma(F)},$$
and we apply it as follows. Fix any $(z,y)\in K\subset B_{\mu_t}\times F$, and let $\sigma=g_{yy}|_{\{z\}\times F}$, a section of $E=(T^*F)^{\otimes 2}$ and apply this estimate to get
\[\begin{split}
\sup_{\{z\}\times F}|\nabla_yg_{yy}|&\leq C[\nabla_y g_{yy}]_{C^\gamma(\{z\}\times F,\omega_F)}\leq C\delta_t^{3+\gamma}[\nabla_y g_{yy}]_{C^\gamma(\{z\}\times F,\ti{\omega}_t)}\\
&\leq C\delta_t^{3+\gamma}\|\omega_t\|_{C^{1,\gamma}(K,\ti{\omega}_t)}\leq C\delta_t^{3+\gamma},
\end{split}\]
using \eqref{bega} and so $\delta_t^{-3}|\nabla_y g_{yy}|\leq C\delta_t^\gamma\to 0$ as claimed, and similarly for the other terms.

Lastly, an argument adapted from \cite{To1} shows that $\omega_0^m$ equals a constant times the Euclidean volume form, and so by standard bootstrapping $\omega_0$ is smooth and Ricci-flat. Thus, we again reach a contradiction to the Liouville Theorem for the complex Monge-Amp\`ere equation in $\mathbb{C}^m$ (see \cite{HT2}). This concludes the sketch of proof of \eqref{bersa} for $k=1$.

\subsection{Setup I.B.2: no fiber space}\label{below}
In this section we discuss the second subcase of the collapsing Setup I.B, which is the following: we have a Calabi-Yau manifold $X^n$ with a cohomology class $[\alpha]\in\de\mathcal{C}_X$ with $\int_X\alpha^n=0$, which we approach with the path $[\alpha_t]=[\alpha]+t[\omega_X]$ (where $\omega_X$ is a K\"ahler metric on $X$), we let $\omega_t$ be the unique Ricci-flat K\"ahler metric in the class $[\alpha_t]$ for $t>0$, but now (unlike in Setup I.B.1) we assume that $[\alpha]$ is not the pullback of a K\"ahler class from the base of a fiber space structure on $X$. We will call this Setup I.B.2.

The Ricci-flat metrics $\omega_t$ are therefore still volume collapsing, but there is no fibration structure so most of the techniques used in Setup I.B.1 are not available anymore. 

The simplest examples are ``irrational'' classes on the boundary of the K\"ahler cone of a torus $X=\mathbb{C}^n/\Lambda$, which are represented by a constant Hermitian semipositive $(1,1)$-form $\beta$ whose kernel foliation has some non-closed leaf (equivalently, $\ker\beta$ viewed as a linear subspace of $\mathbb{C}^n$ is not $\mathbb{Q}$-defined modulo $\Lambda$). In this case, the flat metrics $\omega_t$ on $X$ simply converge smoothly to $\beta$, while the Gromov-Hausdorff limit of $(X,\omega_t)$ is a real torus $T^k$ where $2n-k$ is the real dimension of the largest leaf closure.

However, in general, the behavior can be much worse as shown by Filip and the author \cite{FT,FT2}:

\begin{theorem}\label{dyn}
There are $K3$ surfaces $X$ with classes $[\alpha]\in\de\mathcal{C}_X$ with $\int_X\alpha^2=0$ that satisfy Setup I.B.2 such that the Ricci-flat metrics $\omega_t$ converge weakly as currents to $\eta$, the unique closed positive current in $[\alpha]$, but $\eta$ is not even continuous on the complement of any closed analytic subvariety. In particular, the Ricci-flat metrics $\omega_t$ cannot have a uniform $C^\gamma_{\rm loc}(X\backslash S,\omega_X)$ bound for any $\gamma>0$ and any closed analytic subvariety $S\subset X$.
\end{theorem}

This result in particular implies that such classes $[\alpha]$ do not admit any smooth semipositive representative. There are examples where the $K3$ surface is projective and the class $[\alpha]$ equals $c_1(D)$ for some nef real $\mathbb{R}$-divisor $D$.

The negative result in Theorem \ref{dyn} contrasts with the positive results that we discussed in Setup I.B.1, and show that once the fiber space structure is lost, the Ricci-flat metrics can be rather wildly behaved.

\subsubsection{Holomorphic dynamics} The construction of such $X$ and $[\alpha]$ uses holomorphic dynamics \cite{Can,McM}: these $K3$ surfaces support a chaotic automorphism $T$, in the sense that its topological entropy $h$ is positive, and the class $[\alpha]$ is an eigenclass for $T^*$ acting on $H^{1,1}(X,\mathbb{R})$, with eigenvalue $e^h>1$
$$T^*[\alpha]=e^h[\alpha].$$
Before we explain the idea of how this is used to prove Theorem \ref{dyn}, let us give a few examples of such automorphisms:
\begin{example}\label{exk}
Consider the torus $Y=\mathbb{C}^2/\Lambda$ where $\Lambda=(\mathbb{Z}\oplus i\mathbb{Z})^2$ is the ``square'' lattice, and let $T_Y:Y\to Y$ be the automorphism induced by the linear map $\begin{pmatrix}
2&1\\ 
1&1\end{pmatrix},$ which is known as ``Arnol'd's cat map''. An elementary computation shows that the largest eigenvalue of $T_Y^*$ on $H^{1,1}(Y,\mathbb{R})$ is $e^h=\left(\frac{3+\sqrt{5}}{2}\right)^2$. If we pass to the associated Kummer $K3$ surface $X$, then $T_Y$ induces an automorphism $T$ of $X$ with the same value of topological entropy. The eigenclass $[\alpha]$ with eigenvalue $e^h$ (which is unique up to scaling) admits a smooth semipositive representative.
\end{example}

\begin{example}\label{ex2}
Let $X$ be a $K3$ surface which is a generic complete intersection of degree $(2,2,2)$ in $\mathbb{CP}^1\times\mathbb{CP}^1\times\mathbb{CP}^1$. Then the three projections to $\mathbb{CP}^1\times\mathbb{CP}^1$ exhibit $X$ as a ramified $2:1$ cover, and $X$ is then equipped with three involutions $\sigma_j,j=1,2,3,$ that permute the sheets. The composition $T=\sigma_1\circ\sigma_2\circ\sigma_3$ has $e^h=9+4\sqrt{5}.$ In this case, the eigenclass $[\alpha]$ can be chosen to be $c_1(D)$ for some nef real $\mathbb{R}$-divisor $D$. This construction first appeared in \cite{Maz}.
\end{example}

\begin{example}\label{exm}
By using the Torelli Theorem, McMullen constructed \cite{McM} examples of non-projective $K3$ surfaces $X$ with an automorphism $T$ with $h>0$ which admit a Siegel disc, i.e. an open subset $\Delta\subset X$ biholomorphic to the bidisc in $\mathbb{C}^2$ which is preserved by $T$ and where $T$ is holomorphically conjugate to an irrational rotation $(z_1,z_2)\mapsto (a z_1,bz_2)$ with $|a|=|b|=1$ and $a,b$ and $ab$ are not roots of unity.
\end{example}

If now $T:X\to X$ is a $K3$ automorphism with largest eigenvalue of $T^*$ on $H^{1,1}(X,\mathbb{R})$ larger than $1$, so equal to $e^h,h>0$, then Cantat \cite{Can} shows that there are two nontrivial nef classes $[\eta_{\pm}]\in \de\mathcal{C}_X$ with $\int_X\eta_{\pm}^2=0,\int_X\eta_+\wedge\eta_-=1$ which satisfy
$$T^*[\eta_{\pm}]=e^{\pm h}[\eta_{\pm}],$$
and are not the pullback of a K\"ahler class from the base of an elliptic fibration, so we are in Setup I.B.2. Note that replacing $T$ with $T^{-1}$ interchanges the roles of $[\eta_+]$ and $[\eta_-]$.

Fix smooth closed real $(1,1)$-forms $\alpha_{\pm}$ which represent $[\eta_{\pm}]$. Then Cantat \cite{Can} shows that the classes $[\eta_{\pm}]$ contain unique closed positive currents, which are of the form $\eta_\pm=\alpha_{\pm}+\ddbar\vp_{\pm}\geq 0$ for quasi-psh functions $\vp_{\pm}$ which belong to $C^\gamma(X)$ for some $\gamma>0$ by \cite{DiS}. By uniqueness, they satisfy $T^*\eta_{\pm}=e^{\pm h}\eta_{\pm}$, where the pullback of currents here is defined by $T^*\eta_{\pm}=T^*\alpha_{\pm}+\ddbar(\vp_\pm\circ T)$.

We now take (as we did above) $[\alpha]=[\eta_+]$, and consider as in Setup I.B.2 the Ricci-flat K\"ahler metrics $\omega_t$ on $X$ with cohomology class $[\alpha]+t[\omega_X], t>0$. By weak compactness of closed positive $(1,1)$ currents with bounded cohomology class, from every sequence $t_i\to 0$ one can extract a subsequence (still denoted by $t_i$) such that $\omega_{t_i}$ converge in the weak topology of currents to some closed positive $(1,1)$ current in the class $[\alpha]=[\eta_+]$. But such a current is unique, as mentioned above, and so we see that the whole family $\omega_t$ converges weakly to $\eta_+$ as $t\to 0$.

Suppose that both $\eta_+$ and $\eta_-$ were smooth, so they are closed real semipositive $(1,1)$-forms. Then their wedge product $\mu=\eta_+\wedge\eta_-$ would be a $T$-invariant volume form $\mu$ with total mass $1$, and Cantat \cite{Can} shows that it is also $T$-ergodic. On the other hand, if $\Omega$ is a never-vanishing holomorphic $2$-form on $X$, normalized so that $\mathrm{dVol}=\Omega\wedge\ov{\Omega}$ is a unit-volume volume form, then $\mathrm{dVol}$ is also easily seen to be $T$-invariant as well. We can then write $\mu=f\mathrm{dVol}$ for some smooth nonnegative function $f$ (with integral $1$ with respect to $\mathrm{dVol}$), which is $T$-invariant and hence constant $\mu$-a.e. by ergodicity of $\mu$. But since $\mu=f\mathrm{dVol}$, this means that $f$ is constant on the set where it is strictly positive (which is nonempty), and since $f$ is smooth this implies that $f$ is a constant everywhere, and this constant must be $1$. We thus conclude that $\mu=\mathrm{dVol}$. With small changes to the argument, one obtains the same conclusion also by just assuming that $\eta_\pm$ are both continuous off of some closed analytic subvariety of $X$.

Now, in the Kummer example \ref{exk} there is no contradiction, since indeed here we have $\mu=\mathrm{dVol}$. However we claim that this is not the case in both examples \ref{ex2} and \ref{exm}, which therefore shows that at least one among $\eta_+$ and $\eta_-$ cannot be continuous on the complement of any subvariety, which would prove Theorem \ref{dyn} up to replacing $T$ by $T^{-1}$.

For Example \ref{exm} we argue as follows: it is easy to see in general that for any given K\"ahler metric the iterates
$$\frac{(T^n)^*\omega}{e^{nh}}\to \eta_+,$$
weakly as $n\to+\infty$ (this is because this holds at the level of cohomology classes, and $\eta_+$ is the unique closed positive current it its class). But now on the Siegel disc $\Delta$ we have that $T|_{\Delta}$ is conjugate to a rotation, so on any rotation-invariant compact $K\subset\Delta$ we can choose a $T$-invariant Riemannian metric $g$, and so
$$\sup_K\left|\frac{(T^n)^*\omega}{e^{nh}}\right|_g=e^{-nh}\sup_K|\omega|_g\leq Ce^{-nh},$$
which converges to $0$, thus showing that $\eta_+|_{\Delta}=0$ and so $\mu|_{\Delta}=0$ which is a contradiction to $\mu=\mathrm{dVol}$.

This argument cannot work for Example \ref{ex2} since these $K3$ surfaces are projective and so $T$ cannot have a Siegel disc, as shown in \cite{McM}. We use instead the following rigidity result, proved in the projective case by Cantat-Dupont \cite{CD} and in general by Filip and the author \cite{FT2}:

\begin{theorem}\label{rig}
Let $T:X\to X$ be an automorphism of a $K3$ surface with positive topological entropy, and let $\mu=\eta_+\wedge\eta_-$ and $\mathrm{dVol}=\Omega\wedge\ov{\Omega}$ be the two invariant measures described above. Then the following are equivalent:
\begin{itemize}
\item[(a)] $\mu\ll\mathrm{dVol}$
\item[(b)] $\mu=\mathrm{dVol}$
\item[(c)] $X$ is Kummer and $T$ is induced by an affine linear automorphism of the torus
\end{itemize}
\end{theorem}
To clarify, the measure $\mu$ can be defined in general as the wedge product of $\eta_+$ and $\eta_-$ (not just when these are smooth) because these currents have H\"older continuous potentials and their wedge product was defined by Bedford-Taylor \cite{BT} as a $T$-invariant Radon measure $\mu$ on $X$.

To conclude our argument, it suffices to note that the $K3$ surfaces in Example \ref{ex2} have Picard number $3$ and so cannot be Kummer (which have Picard number at least $16$, coming from the exceptional divisors), and so using Theorem \ref{rig} we obtain a contradiction to $\mu=\mathrm{dVol}$. 

\subsubsection{Other boundary classes}\label{later} Apart from the negative results that we just discussed, essentially nothing is known about the behavior of the Ricci-flat K\"ahler metrics in Setup I.B.2 above, where the limiting nef class $[\alpha]$ has vanishing self-intersection and does not come from the base of a fiber space.
We will now restrict our attention to the case of $K3$ surfaces. The examples above coming from dynamics produce certain nef classes $[\alpha]$ as above, which as we said contain only one closed positive $(1,1)$ current which has a H\"older continuous (but in general not smooth) potential.

This uniqueness property is expected to hold for all nef classes $[\alpha]$ with vanishing self-intersection which don't come from the base of an elliptic fibration (or equivalently, none of the rescaled classes $\lambda[\alpha], \lambda\in\mathbb{R}_{>0},$ belongs to $H^2(X,\mathbb{Q})$), see the ongoing work of Sibony-Verbitsky \cite{SV}.

On the other hand, recall that according to Conjecture \ref{2} we expect that every such class $[\alpha]$ has a closed positive $(1,1)$ current representing it (which is expected to be unique) with $L^\infty$ potential, i.e. of the form $\alpha+\ddbar\vp\geq 0$ where $\alpha$ is a smooth representative of the class and $\vp$ is a bounded quasi-psh function (which can be normalized by $\sup_X\vp=0$). 

We can also write our Ricci-flat metrics as $\omega_t=\alpha+t\omega_X+\ddbar\vp_t,$ $\sup_X\vp_t=0$, for smooth K\"ahler potentials $\vp_t$, and then if Conjecture \ref{2} was settled, then from estimates in pluripotential theory in \cite{FGS} generalizing \cite{BEGZ,DPa,EGZ2,Kol} it would follow that $\|\vp_t\|_{L^\infty(X)}\leq C$ as $t\to 0$ (in all dimensions, not just on $K3$ surfaces).

Furthermore, in the case of $K3$ surfaces, we expect the optimal regularity for the limiting potential $\vp$ to be $C^0$ and no better (in general, despite the fact that in the dynamical examples above the potential was H\"older continuous).

Lastly, one would also like to understand the Gromov-Hausdorff limit of $(X,\omega_t)$ (in particular, show it exists independent of subsequences). In the examples above from dynamics, this limit is in fact just a point, even though the class $[\alpha]$ is not the trivial class. This is easiest to see in Example \ref{ex2} when we take the sequence of Ricci-flat K\"ahler metrics $\frac{(T^n)^*\omega_X}{e^{nh}}$ (which for suitably chosen $\omega_X$ will move towards $[\alpha]=[\eta_+]$ along the same path as $\omega_t$), since these are obviously isometric to $e^{-nh}\omega_X$ and so shrink to a point. 

We then expect the following to hold:
\begin{conjecture}
For every $K3$ surface $X$ and class $[\alpha]\in\de\mathcal{C}_X$ in Setup I.B.2, the Ricci-flat manifolds $(X,\omega_t)$ converge to a point in the Gromov-Hausdorff topology.
\end{conjecture}

\section{Degenerating the complex structure}
In this section we will discuss Setup II from \S \ref{sect2}. Recall that we have now a projective family $\pi:\mathfrak{X}\to\Delta$, smooth over $\Delta^*$, with polarization $\mathfrak{L}$ and with Ricci-flat K\"ahler metrics $\omega_t\in c_1(\mathfrak{L}|_{X_t})$ on $X_t=\pi^{-1}(t)$ for $t\in \Delta^*$. Up to passing to a finite cover, we will assume for simplicity that $K_{X_t}\cong\mathcal{O}_{X_t}$ for all $t\in\Delta^*$.

We can consider the normalization $\nu:\ti{\mathfrak{X}}\to \mathfrak{X}$ and the composition $\ti{\pi}=\pi\circ\nu$, and then if $K_{\ti{\mathfrak{X}}}$ denotes the canonical sheaf then $\ti{\pi}_*K_{\ti{\mathfrak{X}}/\Delta}$ is a line bundle on $\Delta$, which therefore admits a trivializing section, which under the isomorphism $H^0(\Delta, \ti{\pi}_*K_{\ti{\mathfrak{X}}/\Delta})\cong H^0(\ti{\mathfrak{X}},K_{\ti{\mathfrak{X}}})$ corresponds to a holomorphic section $\Omega\in H^0(\ti{\mathfrak{X}},K_{\ti{\mathfrak{X}}})$. Since $\nu$ is an isomorphism away from the central fiber, we can restrict $\Omega$ to $X_t,t\neq 0$, and obtain a trivialization $\Omega_t$ of $K_{X_t}$.

Then the Ricci-flat K\"ahler metric $\omega_t$ on $X_t$ satisfies
\begin{equation}\label{ma3}
\omega_t^n=\left(\frac{\int_{X_t}c_1(\mathfrak{L}|_{X_t})^n}{\int_{X_t}(\sqrt{-1})^{n^2}\Omega_t\wedge\ov{\Omega_t}}\right)(\sqrt{-1})^{n^2}\Omega_t\wedge\ov{\Omega_t}.
\end{equation}

\subsection{Volume}
The volume of $(X_t,\omega_t)$ equals the intersection number $\int_{X_t}c_1(\mathfrak{L}|_{X_t})^n$ which is a constant independent of $t$ (since $\pi$ is flat).

\subsection{Diameter} Unlike in Setup I, where the diameter of $\omega_t$ was always bounded above thanks to Theorem \ref{diam}, the diameter of $(X_t,\omega_t)$ in Setup II may go to infinity. In fact, we have the following result, which combines work of Wang \cite{Wa}, Takayama \cite{Tak} and the author \cite{To5} (when the base of the family is a quasi-projective curve, the stated case with base a disc follows by using also \cite{KNX}), see also the exposition in \cite{Zh3}:

\begin{theorem}\label{taka}
In Setup II, the following are equivalent:
\begin{itemize}
\item $\diam(X_t,\omega_t)\leq C$ as $t\to 0$
\item There is a new family $\mathfrak{X}'$, obtained from $\mathfrak{X}$ by a finite base-change and a birational modification along the central fiber, such that $X'_0$ is a normal Calabi-Yau variety with at worst canonical singularities
\item $0\in\Delta$ is at finite distance from $\Delta^*$ equipped with the Weil-Petersson pseudometric $\omega_{\rm WP}=-\ddbar\log\left(\int_{X_t}(\sqrt{-1})^{n^2}\Omega_t\wedge\ov{\Omega_t}\right)$.
\end{itemize}
\end{theorem}
More precisely, the new family $\mathfrak{X}'$ is obtained from $\mathfrak{X}$ by applying semi-stable reduction followed by a relative Minimal Model Program. 

Since we prefer to work with Ricci-flat metrics with bounded diameter, so that their Gromov-Hausdorff limits will be compact, so we will consider the rescaled metrics
$$\ti{\omega}_t=\frac{\omega_t}{\diam(X_t,\omega_t)^2},$$
which have unit diameter. We therefore see that the volume of $(X_t,\ti{\omega}_t)$ remains bounded away from $0$ if and only if one of the equivalent conditions in Theorem \ref{taka} holds. So, as in Setup I, we have to deal separately with the volume non-collapsed case and with the volume collapsed case.

\subsection{Setup II.A: volume non-collapsed}
Here we assume that we are in Setup II and that any of the equivalent conditions in Theorem \ref{taka} holds, thus the Ricci-flat manifolds $(X_t,\omega_t)$ have fixed volume and diameter bounded uniformly bounded above (and also from below away from $0$ thanks to the Bishop-Gromov volume comparison theorem). Replacing $\mathfrak{X}$ with $\mathfrak{X}'$ (which doesn't change the fibers away from $0$), we may assume that $X_0$ itself is a normal Calabi-Yau variety with at worst canonical singularities.

Thanks to a result of Eyssidieux-Guedj-Zeriahi \cite{EGZ}, which extends the pioneering work of Yau \cite{Ya}, $X_0$ admits a unique ``singular Ricci-flat K\"ahler metric'' $\omega_0$ in $c_1(\mathfrak{L}|_{X_0})$. This can be understood as a weak solution of the corresponding complex Monge-Amp\`ere equation in the sense of pluripotential theory, but all that is important for us is that on the regular part $X_0^{\rm reg}$, $\omega_0$ defines an honest Ricci-flat K\"ahler metric.

Then Rong-Zhang \cite{RZ} prove the following:

\begin{theorem}\label{rz} Suppose in Setup II we assume that $X_0$ is a normal Calabi-Yau variety with at worst canonical singularities. Then $(X_t,\omega_t)$ converge to $(X^{\rm reg}_0,\omega_0)$ locally smoothly on compact subsets of $X^{\rm reg}_0$. Furthermore, the metric completion $(Z,d)$ of $(X^{\rm reg}_0,\omega_0)$ is a compact metric space and $(X_t,\omega_t)\to (Z,d)$ in the Gromov-Hausdorff topology as $t\to 0$.
\end{theorem}
To clarify the statement, locally smooth convergence means that there is a smooth fiber-preserving embedding $F:X^{\rm reg}_0\times\Delta\hookrightarrow\mathfrak{X}$, which is the identity map on $X^{\rm reg}_0\times\{0\}$, such that $(F|_{X^{\rm reg}_0\times\{t\}})^*\omega_t\to \omega_0$ in the smooth topology on compact subsets of $X^{\rm reg}_0$.

The Gromov-Hausdorff limit $Z$ in fact has an algebraic structure, as shown by Donaldson-Sun \cite{DS} (see also \cite[Lemma 2.2]{Zh}):
\begin{theorem}Under the same assumptions as Theorem \ref{rz}, we have that $Z$ is homeomorphic to $X_0$. The homeomorphism identifies the singular set of $Z$ (in the sense of \cite{CC1}) with $X_0\backslash X_0^{\rm reg}$.
\end{theorem}

There are of course many examples of families $\mathfrak{X}$ that satisfy our hypothesis on $X_0$, for example the family of $K3$ surfaces in $\mathbb{CP}^3\times\Delta$ given (in affine coordinates) by
$$x^4+x^2+y^4+y^2+z^4+z^2=t,$$
whose central fiber $X_0$ is a nodal $K3$ surface.

These results give an answer to Question \ref{mainq} in Setup II.A.

\subsection{Setup II.B: volume collapsed}
Here we assume that we are in Setup II and that $\diam(X_t,\omega_t)\to\infty$, i.e. the negation of any of the equivalent conditions in Theorem \ref{taka}. As before, we consider the rescaled Ricci-flat metrics $\ti{\omega}_t=\diam(X_t,\omega_t)^{-2}\omega_t$, so that the Ricci-flat manifolds $(X_t,\ti{\omega}_t)$ have unit diameter and their total volume is collapsing to zero.

An algebro-geometric construction due to Kontsevich-Soibelman \cite{KS2}, refined in \cite{MN,NX}, associates to our family $\mathfrak{X}$ a connected simplicial complex $Sk(\mathfrak{X})$, the essential skeleton, which is roughly speaking the dual intersection complex of the central fiber $X'_0$ of the family $\mathfrak{X}'$ in Theorem \ref{taka}. The real dimension of the essential skeleton satisfies $0\leq\dim_{\mathbb{R}}Sk(\mathfrak{X})\leq n$, and we have that $\dim_{\mathbb{R}}Sk(\mathfrak{X})=0$ if and only if $X'_0$ is a normal Calabi-Yau variety with at worst canonical singularities, and that $\dim_{\mathbb{R}}Sk(\mathfrak{X})=n$ if and only if $\mathfrak{X}$ is a {\em maximal degeneration} or {\em large complex structure limit} in the sense of mirror symmetry, see e.g. \cite{GW,Mor}.

The simplest example of a large complex structure limit family is the family of elliptic curves in $\mathbb{CP}^2\times\Delta$ given (in affine coordinates) by
$$y^2=x^3+x^2+t,$$
whose central fiber $X_0$ is a nodal elliptic curve. The essential skeleton of this family is $S^1$. In this case we can see the behavior of $\ti{\omega}_t$ explicitly: for $t\neq 0$ we can write $X_t=\mathbb{C}/\Lambda_t$ where $\Lambda_t$ is spanned over $\mathbb{Z}$ by $1$ and $\tau(t)$, and an explicit calculation shows that as $t\in\mathbb{R}$ approaches $0$ then $\tau(t)$ is approximately equal to $i|\log t|$. This implies that the unit-diameter rescaled flat metrics $\ti{\omega}_t$ come from the Euclidean metric with a fundamental domain which is approximately a rectangle with sides $1$ and $\tau(t)^{-1}$, and so their Gromov-Hausdorff limit is the unit-diameter $S^1$.

Another standard example of a large complex structure limit family is the one in $\mathbb{CP}^{n+1}\times\Delta$ given (in homogeneous coordinates) by
$$z_0^{n+2}+\cdots+z_{n+1}^{n+2}+\frac{1}{t}z_0\cdots z_{n+1}=0,$$
with smooth Calabi-Yau $n$-folds $X_t$ for $t\neq 0$ and with $X_0=\{z_0\cdots z_{n+1}=0\}$ the union of the $n+2$ coordinate hyperplanes. Note that the dual intersection complex of $X_0$ is the boundary of the standard $(n+1)$-simplex in $\mathbb{R}^{n+2}$, and so homeomorphic to $S^n$. In fact, this is the essential skeleton of $\mathfrak{X}$, i.e. $Sk(\mathfrak{X})\cong S^n$.

It is conjectured (see \cite{KLSV,KX}) that for large complex structure limits where the Calabi-Yau fibers $X_t$ are simply connected with $H^i(X_t,\mathcal{O}_{X_t})=0, 0<i<n,$ the skeleton $Sk(\mathfrak{X})$ should be homeomorphic to $S^n$, while when $X_t$ is hyperk\"ahler it should be homeomorphic to $\mathbb{CP}^{\frac{n}{2}}$.

\subsection{Large complex structure limits}
The guiding conjecture about large complex structure limits of Calabi-Yau manifolds is due to Strominger-Yau-Zaslow \cite{SYZ} (SYZ):
\begin{conjecture}\label{syz}
Given a large complex structure limit family $\pi:\mathfrak{X}\to\Delta$ of Calabi-Yau $n$-folds, for $t\neq 0$ sufficiently small the manifold $X_t$ admits a Special Lagrangian $T^n$-fibration over a half-dimensional base $B$, with singular fibers lying over a subset $S\subset B$ of codimension $2$.
\end{conjecture}
By taking the fiberwise dual tori one obtains dual $T^n$-fibrations over $B\backslash S$, which should be compactifiable to the ``mirror family'' $\check{\mathfrak{X}}$ of $\mathfrak{X}$ (which satisfy certain properties that we will not get into, but see e.g. \cite{GHJ} for an introduction). The base $B$ is expected to be homeomorphic to $S^n$ if the Calabi-Yau fibers $X_t$ are simply connected with $H^i(X_t,\mathcal{O}_{X_t})=0, 0<i<n,$ and to be homeomorphic to $\mathbb{CP}^{\frac{n}{2}}$ if $X_t$ is hyperk\"ahler (there is also strictly speaking the elementary case when $X_t$ are tori, where $B$ is also $T^n$).

Conjecture \ref{syz} remains open, but see \cite{Li4} for a recent result in this direction, and also the Gross-Siebert program \cite{GS} for an algebro-geometric approach to the construction of $\check{\mathfrak{X}}$.

Around 2000, Kontsevich-Soibelman \cite{KS,KS2}, Gross-Wilson \cite{GW} and Todorov \cite{Man} proposed that the base $B$ of the conjectural SYZ fibration arises as the collapsed Gromov-Hausdorff limit of Ricci-flat metrics on $X_t$, exactly in the setup II.B that we have been considering. More precisely, they conjectured:

\begin{conjecture}\label{ks}
Given a large complex structure limit family $\pi:\mathfrak{X}\to\Delta$ of Calabi-Yau $n$-folds, equip $X_t,t\neq 0,$ with the Ricci-flat K\"ahler metrics $\omega_t\in c_1(\mathfrak{L}|_{X_t})$ and consider the unit-diameter rescaled metrics $\ti{\omega}_t=\diam(X_t,\omega_t)^{-2}\omega_t$. Then as $t\to 0$, $(X_t,\ti{\omega}_t)\to (Z,d)$ in the Gromov-Hausdorff topology, where $(Z,d)$ is a compact metric space with an open dense subset $Z_0$ over which $d|_{Z_0}$ is induced by a Riemannian metric $g_0$, and the real Hausdorff codimension of $Z\backslash Z_0$ is at least $2$. Furthermore, $Z_0$ carries an integral affine structure and in local affine coordinates $g_0=\nabla^2 F$ where $F$ is a smooth convex function that satisfies the real Monge-Amp\`ere equation $\det(\nabla^2 F)=c\in\mathbb{R}_{>0}$.
\end{conjecture}
In other words, $(Z_0,g_0)$ is a Monge-Amp\`ere manifold in the terminology of Cheng-Yau \cite{CY}. Of course, $Z$ is expected to be equal to the base $B$ of the conjectural SYZ fibration, and therefore it is also expected to be topologically $S^n$ or $\mathbb{CP}^{\frac{n}{2}}$ (as above), and so it is expected to be homeomorphic to $Sk(\mathfrak{X})$. Furthermore, passing to the dual affine structure, and applying Legendre transform to $F$, one obtains a dual Monge-Amp\`ere manifold $(Z_0,\check{g}_0)$. From this dual Monge-Amp\`ere manifold it is elementary to construct a family of smooth $T^n$-fibration over $Z_0$ with fiber size depending on $t$, which for $t$ small is expected to compactify across $Z\backslash Z_0$ to the mirror family $\check{X}_t$, see \cite{KS2}.

Going back to Conjecture \ref{ks}, which fits precisely in Question \ref{mainq}, it is not hard to see that it holds when $X_t$ are tori (see \cite{Od}). For other Calabi-Yau manifolds, it was proved by Gross-Wilson \cite{GW} for those large complex structure limits of $K3$ surfaces which arise from hyperk\"ahler rotation from elliptically fibered $K3$ surfaces (as in Setup I.B.1) with only $I_1$ singular fibers. Later the author with Gross and Zhang \cite{GTZ,GTZ2} gave a new proof of this result, which allows for arbitrary singular fibers (see also \cite{CVZ} for a very recent extension of the Gross-Wilson gluing method), and the author with Zhang \cite{TZ2} proved it for some families of hyperk\"ahler manifolds which arise from hyperk\"ahler rotation from hyperk\"ahler manifolds with a holomorphic Lagrangian fibration; these results were described in the section about Setup I.B.1. In all of these cases, we also proved in \cite{GTZ2,TZ2} that the limiting manifold $(Z_0,g_0)$ is {\em special K\"ahler} in the sense of \cite{Fr, Str}. Using these results, Odaka-Oshima \cite{OO} proved Conjecture \ref{ks} in general for $K3$ surfaces. Lastly, the first step towards a Gross-Wilson type gluing construction for large complex structure limits of Calabi-Yau $3$-folds was taken in \cite{Li3}.

In a related direction, Boucksom-Jonsson \cite{BJ} have investigated the fiberwise Ricci-flat volume forms $\omega_t^n$ on $X_t$, and they proved that $\omega_t^n$ converges in the weak topology to a Lebesgue-type measure $\mu$ on the essential skeleton $Sk(\mathfrak{X})$. To make sense of this convergence, they construct a {\em hybrid space} $\mathfrak{X}^{\rm hyb}$ over $\Delta$ which coincides with $\mathfrak{X}$ over $\Delta^*$ but over $0$ one has the Berkovich analytification $\mathfrak{X}^{\rm an}$ associated to the basechange of $\mathfrak{X}$ over $\Delta^*$ to the non-Archimedean field $\mathbb{C}((t))$, inside which $Sk(\mathfrak{X})$ naturally sits. On the other hand, Boucksom-Favre-Jonsson \cite{BFJ} have managed to solve the non-Archimedean analog of the complex Monge-Amp\`ere equation on $Sk(\mathfrak{X})$ with volume form $\mu$ (when $\mathfrak{X}$ is defined over an algebraic curve, with \cite{BGJK} later removing this assumption), and it is natural to ask about the relation between their solution and Conjecture \ref{ks}, and whether this circle of ideas can be used to attack it.

\subsection{Other collapsing limits}
Lastly, we consider the case of a family $\mathfrak{X}$ of Calabi-Yau $n$-folds as in Setup II.B such that $0<\dim_{\mathbb{R}}Sk(\mathfrak{X})<n$. In this case the following is expected:

\begin{conjecture}\label{skel}
Given a family $\mathfrak{X}\subset\mathbb{CP}^N\times\Delta$ of Calabi-Yau $n$-folds with $0<\dim_{\mathbb{R}}Sk(\mathfrak{X})<n$, equip $X_t,t\neq 0,$ with the Ricci-flat K\"ahler metrics $\omega_t\in c_1(\mathfrak{L}|_{X_t})$ and consider the unit-diameter rescaled metrics $\ti{\omega}_t=\diam(X_t,\omega_t)^{-2}\omega_t$. Then as $t\to 0$, $(X_t,\ti{\omega}_t)\to (Z,d)$ in the Gromov-Hausdorff topology, where $(Z,d)$ is a compact metric space homeomorphic to the essential skeleton $Sk(\mathfrak{X})$.
\end{conjecture}

The first nontrivial examples are families of $K3$ surfaces with $\dim_{\mathbb{R}}Sk(\mathfrak{X})=1$, for example the family $\mathfrak{X}$ in $\mathbb{CP}^{3}\times\Delta$
\begin{equation}\label{k3a}
P_1(z)P_2(z)+tP(z)=0,
\end{equation}
where $P_1,P_2$ are generic quadratic polynomials and $P$ is a generic quartic polynomial. The central fiber $X_0$ has two irreducible components with nontrivial intersection, so its dual intersection complex is homeomorphic to the interval $[0,1]$, and so is $Sk(\mathfrak{X})$.

The only examples where Conjecture \ref{skel} is known are the families considered in \cite{SZ}, for which the Ricci-flat metrics $\ti{\omega}_t$ can be obtained via a gluing construction. These include the $K3$ families in \eqref{k3a}, where collapsing to a $1$-dimensional interval is proved. Note that there are earlier gluing constructions of collapsing hyperk\"ahler metrics on $K3$ surfaces in \cite{HSVZ} which converge to $[0,1]$ but it is not clear whether some of them come from projective families as in Conjecture \ref{skel}. There are also other examples in \cite{Fos} of hyperk\"ahler metrics on $K3$ surfaces which collapse to $T^3/\mathbb{Z}_2$, and will not arise from Conjecture \ref{skel}. We refer the reader to the recent survey \cite{Su}.

Lastly, we mention one more problem by Kontsevich-Soibelman \cite{KS2} which remains open:
\begin{conjecture}Suppose that $\mathfrak{X}$ is a Calabi-Yau family as in Setup II, which does not satisfy the equivalent conditions in Theorem \ref{taka}. Then there are constants $A,C>0$ such that the Ricci-flat K\"ahler metrics $\omega_t\in c_1(\mathfrak{L}|_{X_t})$ satisfy
\begin{equation}\label{diamblow}
\diam(X_t,\omega_t)\geq C^{-1}(-\log|t|)^A,
\end{equation}
as $t\to 0$.
\end{conjecture}
Of course, by Theorem \ref{taka} we know that $\diam(X_t,\omega_t)\to+\infty$, so the conjecture asks for an effective lower bound. Note that in \cite{KS2} this is only conjectured for large complex structure limit families, but it seems natural to extend this as above. Note also that \eqref{diamblow} holds for the metrics constructed in \cite{SZ} by gluing methods.


\begin{thebibliography}{99}
\bibitem{Be} A. Beauville, {\em Vari\'et\'es K\"ahleriennes dont la premi\`ere classe de Chern est nulle}, J. Differential Geom. {\bf 18} (1983), no. 4, 755--782.
\bibitem{BT} E. Bedford, B.A. Taylor, {\em The Dirichlet problem for a complex Monge-Amp\`ere equation}, Invent. Math. {\bf 37} (1976), no. 1, 1--44.
\bibitem{Bog} F.A. Bogomolov, {\em K\"ahler manifolds with trivial canonical class}, Izv. Akad. Nauk SSSR Ser. Mat. {\bf 38} (1974), 11--21; English translation in Math. USSR Izv. {\bf 8} (1974), no. 1, 9--20.
\bibitem{Bog2} F.A. Bogomolov, {\em The decomposition of K\"ahler manifolds with a trivial canonical class}, Mat. Sb. (N.S.) {\bf 93(135)} (1974), 573--575; English translation in Math. USSR Sbornik {\bf 22} (1974), no. 4, 580--583.
\bibitem{BEGZ} S. Boucksom, P. Eyssidieux, V. Guedj, A. Zeriahi, \emph{ Monge-Amp\`ere equations in big cohomology classes}, Acta Math. {\bf 205} (2010), no. 2, 199--262.
\bibitem{BFJ} S. Boucksom, C. Favre, M. Jonsson, {\em Solution to a non-Archimedean Monge-Amp\`ere equation}, J. Amer. Math. Soc. {\bf 28} (2015), no. 3, 617--667.
\bibitem{BJ} S. Boucksom, M. Jonsson, {\em Tropical and non-Archimedean limits of degenerating families of volume forms}, J. \'Ec. polytech. Math. {\bf 4} (2017), 87--139.
\bibitem{BGJK} J.I. Burgos Gil, W. Gubler, P. Jell, K. K\"unnemann, F. Martin, {\em Differentiability of non-archimedean volumes and non-archimedean Monge-Amp\`ere equations (with an appendix by Robert Lazarsfeld)}, Algebr. Geom. {\bf 7} (2020), no. 2, 113--152.
\bibitem{Cal} E. Calabi, {\em On K\"ahler manifolds with vanishing canonical class}, in {\em Algebraic geometry and topology. A symposium in honor of S. Lefschetz},  pp. 78--89. Princeton University Press, Princeton, N. J., 1957.
\bibitem{Can} S. Cantat, {\em Dynamique des automorphismes des surfaces $K3$}, Acta Math. {\bf 187} (2001), no. 1, 1--57.
\bibitem{CD} S. Cantat, C. Dupont, {\em Automorphisms of surfaces: Kummer rigidity and measure of maximal entropy}, to appear in J. Eur. Math. Soc. (JEMS).
\bibitem{CC1}  J. Cheeger, T.H. Colding, {\em On the structure of spaces with Ricci curvature bounded below. I}, J. Differential Geom. {\bf 46}  (1997),  no. 3, 406--480.
\bibitem{CC2}  J. Cheeger, T.H. Colding, {\em On the structure of spaces with Ricci curvature bounded below. II}, J. Differential Geom. {\bf 54} (2000), no. 1, 13--35.
\bibitem{CC3}  J. Cheeger, T.H. Colding, {\em  On the structure of spaces with Ricci curvature bounded below. III}, J. Differential Geom. {\bf 54} (2000), no. 1, 37--74.
\bibitem{CVZ} G. Chen, J. Viaclovsky, R. Zhang, {\em Collapsing Ricci-flat metrics on elliptic K3 surfaces}, preprint, arXiv:1910.11321.
\bibitem{CY}  S.Y. Cheng, S.T. Yau,  {\em  The real Monge-Amp\`{e}re equation and affine flat structures}, Chern S.S., ed. Proc. 1980 Beijing Sympos. Diff. Geom. Diff. Eq., Vol. 1 (1982), 339--370.
\bibitem{Ch} I. Chiose, {\em The K\"ahler rank of compact complex manifolds}, J. Geom. Anal. {\bf 26} (2016), no. 1, 603--615.
\bibitem{CT} T. Collins, V. Tosatti, {\em K\"ahler currents and null loci}, Invent. Math. {\bf 202} (2015), no.3, 1167--1198.
\bibitem{DJZ} V. Datar, A. Jacob, Y. Zhang, {\em Adiabatic limits of anti-self-dual connections on collapsed $K3$ surfaces}, to appear in J. Differential Geom.
\bibitem{DPa} J.-P. Demailly, N. Pali, \emph{Degenerate complex Monge-Amp\`ere equations over compact K\"ahler manifolds},   Internat. J. Math. {\bf 21} (2010), no. 3, 357--405.
\bibitem{DP} J.-P. Demailly, M. P\u{a}un, \emph{Numerical characterization of the K\"ahler cone of a compact K\"ahler manifold}, Ann. of Math., {\bf 159} (2004), no. 3, 1247--1274.
\bibitem{DPS2} J.-P. Demailly, T. Peternell, M. Schneider, {\em K\"ahler manifolds with numerically effective Ricci class}, Compositio Math. {\bf 89} (1993), no. 2, 217--240.
\bibitem{DPS} J.-P. Demailly, T. Peternell, M. Schneider, {\em Compact complex manifolds with numerically effective tangent bundles}, J. Algebraic Geom. {\bf 3} (1994), no. 2, 295--345.
\bibitem{DiS} T.-C. Dinh, N. Sibony, {\em Green currents for holomorphic automorphisms of compact K\"ahler manifolds}, J. Amer. Math. Soc. {\bf 18} (2005), no. 2, 291--312.
\bibitem{DS} S.K. Donaldson, S. Sun, {\em Gromov-Hausdorff limits of K\"ahler manifolds and algebraic geometry}, Acta Math. {\bf 213} (2014), no. 1, 63--106.
\bibitem{DS2} S.K. Donaldson, S. Sun, {\em Gromov-Hausdorff limits of K\"ahler manifolds and algebraic geometry, II},  J. Differential Geom. {\bf 107} (2017), no. 2, 327--371.
\bibitem{ELMNP} L. Ein, R. Lazarsfeld, M. Musta\c{t}\u{a}, M. Nakamaye, M. Popa, {\em Restricted volumes and base loci of linear series}, Amer. J. Math. {\bf 131} (2009), no. 3, 607--651.
\bibitem{EGZ} P. Eyssidieux, V. Guedj, A. Zeriahi, {\em Singular K\"ahler-Einstein metrics}, J. Amer. Math. Soc. {\bf 22} (2009), no. 3, 607--639.
\bibitem{EGZ2} P. Eyssidieux, V. Guedj, A. Zeriahi, \emph{A priori $L^\infty$-estimates for degenerate complex Monge-Amp\`ere equations}, Int. Math. Res. Not. {\bf 2008}, Art. ID rnn 070, 8 pp.
\bibitem{Fil} S. Filip, {\em An introduction to $K3$ surfaces and their dynamics}, lecture notes.
\bibitem{FT} S. Filip, V. Tosatti, {\em Smooth and rough positive currents}, Ann. Inst. Fourier (Grenoble) {\bf 68} (2018), no.7, 2981--2999.
\bibitem{FT2} S. Filip, V. Tosatti, {\em Kummer rigidity for $K3$ surface automorphisms via Ricci-flat metrics}, to appear in Amer. J. Math.
\bibitem{Fi} J. Fine, {\em Fibrations with constant scalar curvature K\"ahler metrics and the CM-line bundle}, Math. Res. Lett. {\bf 14} (2007), no. 2, 239--247.
\bibitem{FG} W. Fischer, H. Grauert, {\em Lokal-triviale Familien kompakter komplexer Mannigfaltigkeiten}, Nachr. Akad. Wiss. G\"ottingen Math.-Phys. Kl. II (1965), 89--94.
\bibitem{Fos} L. Foscolo, {\em ALF gravitational instantons and collapsing Ricci-flat metrics on the $K3$ surface}, J. Differential Geom. {\bf 112} (2019), no. 1, 79--120.
\bibitem{Fr}   D. Freed,  {\em Special K\"ahler Manifolds},  Comm. Math. Phys. {\bf 203} (1999), no. 1, 31--52.
\bibitem{FGS} X. Fu, B. Guo, J. Song, {\em Geometric estimates for complex Monge-Amp\`ere equations}, to appear in J. Reine Angew. Math.
\bibitem{Fuj} A. Fujiki, {\em On automorphism groups of compact K\"ahler manifolds}, Inv. Math. {\bf 44} (1978), 225--258.
\bibitem{GSVY} B. Greene, A. Shapere, C. Vafa, S.-T. Yau, {\em Stringy cosmic strings and noncompact Calabi-Yau manifolds}, Nuclear Phys. B {\bf 337} (1990), no. 1, 1--36.
\bibitem{GH} P. Griffiths, J. Harris, {\em Principles of algebraic geometry}, John Wiley \& Sons, Inc., New York, 1994.
\bibitem{gr2} M. Gromov, {\em Metric structures for Riemannian and non-Riemannian spaces}, Birkh\"auser, Boston, 2007.
\bibitem{GHJ} M. Gross, D. Huybrechts, D. Joyce, \emph{Calabi-Yau manifolds and related geometries,} Springer-Verlag 2003.
\bibitem{GS} M. Gross, B. Siebert, {\em Intrinsic Mirror Symmetry}, preprint, arXiv:1909.07649.
\bibitem{GTZ} M. Gross, V. Tosatti, Y. Zhang, {\em Collapsing of abelian fibred Calabi-Yau manifolds}, Duke Math. J. {\bf 162} (2013), no. 3, 517--551.
\bibitem{GTZ2} M. Gross, V. Tosatti, Y. Zhang,  {\em Gromov-Hausdorff collapsing of Calabi-Yau manifolds},  Comm. Anal. Geom. {\bf 24} (2016), no. 1, 93--113.
\bibitem{GTZ3} M. Gross, V. Tosatti, Y. Zhang,  {\em Geometry of twisted K\"ahler-Einstein metrics and collapsing}, preprint, arXiv:1911.07315.
\bibitem{GW} M. Gross, P.M.H. Wilson, \emph{Large complex structure limits of $K3$ surfaces}, J. Differ. Geom. {\bf 55} (2000), no. 3, 475--546.
\bibitem{He} H.-J. Hein, {\em Gravitational instantons from rational elliptic surfaces}, J. Amer. Math. Soc. {\bf 25} (2012), no. 2, 355--393.
\bibitem{He2} H.-J. Hein, {\em A Liouville theorem for the complex Monge-Amp\`ere equation on product manifolds}, Comm. Pure Appl. Math. {\bf 72} (2019), no. 1, 122--135.
\bibitem{HS} H.-J. Hein, S. Sun, {\em Calabi-Yau manifolds with isolated conical singularities}, Publ. Math. Inst. Hautes \'Etudes Sci. {\bf 126} (2017), 73--130.
\bibitem{HSVZ} H.-J. Hein, S. Sun, J. Viaclovsky, R. Zhang, {\em Nilpotent structures and collapsing Ricci-flat metrics on $K3$ surfaces}, preprint, arXiv:1807.09367.
\bibitem{HT} H.-J. Hein, V. Tosatti, {\em Remarks on the collapsing of torus fibered Calabi-Yau manifolds},  Bull. Lond. Math. Soc. {\bf 47} (2015), no. 6, 1021--1027.
\bibitem{HT2} H.-J. Hein, V. Tosatti, {\em Higher-order estimates for collapsing Calabi-Yau metrics}, preprint, arXiv:1803.06697.
\bibitem{Hor} A. H\"oring, {\em Adjoint $(1,1)$-classes on threefolds}, preprint, arXiv:1807.08442.
\bibitem{Hw} J.-M. Hwang, {\em  Base manifolds for fibrations of projective irreducible symplectic manifolds}, Invent. Math. {\bf 174} (2008), no. 3, 625--644.
\bibitem{JS} W. Jian, Y. Shi, \emph{A ``boundedness implies convergence'' principle and its applications to collapsing estimates in K\"ahler geometry}, preprint, arXiv:1904.11261.
\bibitem{Kaw} Y. Kawamata, {\em On the length of an extremal rational curve}, Invent. Math. {\bf 105} (1991), no. 3, 609--611.
\bibitem{KMM} Y. Kawamata, K. Matsuda, K. Matsuki, \emph{Introduction to the minimal model problem}, in  \emph{Algebraic geometry, Sendai, 1985}, 283--360, Adv. Stud. Pure Math., 10, North-Holland, Amsterdam, 1987.
\bibitem{KoT} R. Kobayashi, A.N. Todorov, {\em Polarized period map for generalized $K3$ surfaces and the moduli of Einstein metrics}, Tohoku Math. J. (2) {\bf 39} (1987), no. 3, 341--363.
\bibitem{KLSV} J. Koll\'ar, R. Laza, G. Sacc\`a, C. Voisin, {\em Remarks on degenerations of hyper-K\"ahler manifolds}, Ann. Inst. Fourier (Grenoble) {\bf 68} (2018), no. 7, 2837--2882.
\bibitem{KNX} J. Koll\'ar, J. Nicaise, C. Xu, {\em Semi-stable extensions over $1$-dimensional bases}, Acta Math. Sin. (Engl. Ser.) {\bf 34} (2018), no. 1, 103--113.
\bibitem{KX} J. Koll\'ar, C. Xu, {\em The dual complex of Calabi-Yau pairs}, Invent. Math. {\bf 205} (2016), no. 3, 527--557.
\bibitem{Kol} S. Ko\l odziej, {\em The complex Monge-Amp\`ere equation}, Acta Math. {\bf 180} (1998), no. 1, 69--117.
\bibitem{KS} M. Kontsevich, Y. Soibelman, {\em  Homological mirror symmetry and torus fibrations}, in {\em Symplectic geometry and mirror symmetry,} 203--263, World Sci. Publishing 2001.
\bibitem{KS2} M. Kontsevich, Y. Soibelman, {\em  Affine Structures and Non-Archimedean Analytic Spaces}, in \it The Unity of Mathematics, \rm  Progress in Mathematics Volume 244, Springer,  (2006),  321--385.
\bibitem{LOP} V. Lazi\'c, K. Oguiso, T. Peternell, {\em Nef line bundles on Calabi-Yau threefolds, I}, to appear in Int. Math. Res. Not. IMRN
\bibitem{LS} C. LeBrun, M. Singer, {\em A Kummer-type construction of self-dual 4-manifolds}, Math. Ann. {\bf 300} (1994), no. 1, 165--180.
\bibitem{Li} Y. Li, {\em On collapsing Calabi-Yau fibrations}, to appear in J. Differential Geom.
\bibitem{Li2} Y. Li, {\em A gluing construction of collapsing Calabi-Yau metrics on $K3$ fibred $3$-folds}, Geom. Funct. Anal. {\bf 29} (2019), no. 4, 1002--1047.
\bibitem{Li3} Y. Li, {\em SYZ geometry for Calabi-Yau 3-folds: Taub-NUT and Ooguri-Vafa type metrics}, preprint, arXiv:1902.08770.
\bibitem{Li4} Y. Li, {\em SYZ conjecture for Calabi-Yau hypersurfaces in the Fermat family}, preprint, arXiv:1912.02360.
\bibitem{LLZ} C. Li, J. Li, X. Zhang, \emph{A mean value formula and a Liouville theorem for the complex Monge-Amp\`ere equation}, to appear in Int. Math. Res. Not. IMRN.
\bibitem{Lie} D.I. Lieberman, {\em Compactness of the Chow scheme: applications to automorphisms and deformations of K\"ahler manifolds}, in {\em Fonctions de plusieurs variables complexes, III (S\'em. Fran\c{c}ois Norguet, 1975–1977)}, pp. 140--186, Lecture Notes in Math., 670, Springer, Berlin, 1978.
\bibitem{Man} Y.I. Manin, {\em  Moduli, motives, mirrors}, in {\em European Congress of Mathematics, Vol. I (Barcelona, 2000)}, 53--73, Progr. Math., 201, Birkh\"auser, Basel, 2001.
\bibitem{Mat} Y. Matsushima, {\em Holomorphic vector fields and the first Chern class of a Hodge manifold}, J. Differential Geom. {\bf 3} (1969), 477--480.
\bibitem{Maz} B. Mazur, {\em The topology of rational points}, Experiment. Math. {\bf 1} (1992), no. 1, 35--45.
\bibitem{McM} C.T. McMullen, {\em Dynamics on $K3$ surfaces: Salem numbers and Siegel disks}, J. Reine Angew. Math. {\bf 545} (2002), 201--233.
\bibitem{Mor} D. Morrison, {\em Compactifications of moduli spaces inspired by mirror symmetry}, in \it Journ\'{e}es de g\'{e}om\'{e}trie Alg\'{e}brique d'Orsay, \rm Astrisque 218 (1993), 243--271.
\bibitem{MN} M. Musta\c{t}\u{a}, J. Nicaise, {\em Weight functions on non-Archimedean analytic spaces and the Kontsevich-Soibelman skeleton}, Algebr. Geom. {\bf 2} (2015), no. 3, 365--404.
\bibitem{Nak} M. Nakamaye, \emph{Stable base loci of linear series}, Math. Ann., {\bf 318} (2000), no. 4, 837--847.
\bibitem{NX} J. Nicaise, C. Xu, {\em The essential skeleton of a degeneration of algebraic varieties}, Amer. J. Math. {\bf 138} (2016), no. 6, 1645--1667.
\bibitem{Od} Y. Odaka, {\em Tropical geometry compactifications of moduli, II: $A_g$ case and holomorphic limits}, Int. Math. Res. Not. IMRN 2019, no. 21, 6614--6660.
\bibitem{OO} Y. Odaka, Y. Oshima, {\em Collapsing $K3$ surfaces, tropical geometry and moduli compactifications of Satake, Morgan-Shalen type}, preprint, arXiv:1810.07685.
\bibitem{RZ}  X. Rong, Y. Zhang,  {\em   Continuity of Extremal Transitions and Flops for Calabi-Yau Manifolds,}  Appendix B by Mark Gross,  J. Differential Geom. 89  (2011), no. 2,  233--269.
\bibitem{Sat} I. Satake, {\em Algebraic structures of symmetric domains}, Iwanami Shoten, Tokyo; Princeton University Press, Princeton, N.J., 1980.
\bibitem{SY} R. Schoen, S.-T. Yau, {\em Lectures on differential geometry}, International Press, Cambridge, MA, 1994.
\bibitem{SV} N. Sibony, M. Verbitsky, in preparation, \url{http://verbit.ru/MATH/TALKS/Rigid-currents-NYU-2019.pdf}
\bibitem{LSi} L. Simon, {\em Schauder estimates by scaling}, Calc. Var. Partial Differential Equations {\bf 5} (1997), no. 5, 391--407.
\bibitem{So} J. Song, {\em Riemannian geometry of K\"ahler-Einstein currents}, preprint, arXiv:1404.0445.
\bibitem{ST0} J. Song, G. Tian, {\em The K\"ahler-Ricci flow on surfaces of positive Kodaira dimension}, Invent. Math. {\bf 170} (2007), no. 3, 609--653.
\bibitem{ST} J. Song, G. Tian, {\em Canonical measures and K\"ahler-Ricci flow}, J. Amer. Math. Soc. {\bf 25} (2012), no. 2, 303--353.
\bibitem{STZ} J. Song, G. Tian, Z. Zhang, {\em Collapsing behavior of Ricci-flat K\"ahler metrics and long time solutions of the K\"ahler-Ricci flow}, preprint, arXiv:1904.08345.
\bibitem{Str} A. Strominger, {\em Special geometry}, Comm. Math. Phys. {\bf 133} (1990), no. 1, 163--180.
\bibitem{SYZ} A. Strominger, S.-T. Yau, E. Zaslow, {\em Mirror symmetry is $T$-duality}, Nuclear Phys. B {\bf 479} (1996), no. 1-2, 243--259.
\bibitem{Su} S. Sun, {\em Collapsing of Calabi-Yau metrics and degeneration of complex structures}, preprint.
\bibitem{SZ} S. Sun, R. Zhang, {\em Complex structure degenerations and collapsing of Calabi-Yau metrics}, preprint, arXiv:1906.03368.
\bibitem{Tak} S. Takayama, {\em On moderate degenerations of polarized Ricci-flat K\"ahler manifolds}, J. Math. Sci. Univ. Tokyo {\bf 22} (2015), no. 1, 469--489.
\bibitem{To0} V. Tosatti, {\em Limits of Calabi-Yau metrics when the K\"ahler class degenerates}, J. Eur. Math. Soc. (JEMS) {\bf 11} (2009), no. 4, 755--776.
\bibitem{To1} V. Tosatti, {\em  Adiabatic limits of Ricci-flat K\"ahler metrics}, J. Differential Geom. {\bf 84} (2010), no. 2, 427--453.
\bibitem{To2} V. Tosatti, {\em Degenerations of Calabi-Yau metrics}, in {\em Geometry and Physics in Cracow,} Acta Phys. Polon. B Proc. Suppl. {\bf 4} (2011), no. 3, 495--505.
\bibitem{To3} V. Tosatti, {\em Calabi-Yau manifolds and their degenerations}, Ann. N.Y. Acad. Sci. {\bf 1260} (2012), 8--13.
\bibitem{To5} V. Tosatti, {\em Families of Calabi-Yau manifolds and canonical singularities}, Int. Math. Res. Not. IMRN 2015, no. 20, 10586--10594.
\bibitem{To4} V. Tosatti, {\em KAWA lecture notes on the K\"ahler-Ricci flow}, Ann. Fac. Sci. Toulouse Math. {\bf 27} (2018), no. 2, 285--376.
\bibitem{ToS} V. Tosatti, {\em Ricci-flat metrics and dynamics on $K3$ surfaces}, preprint.
\bibitem{TWY} V. Tosatti, B. Weinkove, X. Yang, {\em The K\"ahler-Ricci flow, Ricci-flat metrics and collapsing limits}, Amer. J. Math. {\bf 140} (2018), no. 3, 653--698.
\bibitem{TZ0} V. Tosatti, Y. Zhang, {\em Triviality of fibered Calabi-Yau manifolds without singular fibers}, Math. Res. Lett. {\bf 21} (2014), no. 4, 905--918.
\bibitem{TZ} V. Tosatti, Y. Zhang, {\em Infinite time singularities of the K\"ahler-Ricci flow}, Geom. Topol. {\bf 19} (2015), no. 5, 2925--2948.
\bibitem{TZ2} V. Tosatti, Y. Zhang, {\em Collapsing hyperk\"ahler manifolds}, to appear in Ann. Sci. \'Ec. Norm. Sup\'er.
\bibitem{Tot} B. Totaro, {\em Algebraic surfaces and hyperbolic geometry}, in {\em Current developments in algebraic geometry}, 405--426, Math. Sci. Res. Inst. Publ., 59, Cambridge Univ. Press, Cambridge, 2012.
\bibitem{Wa} C.-L. Wang, {\em On the incompleteness of the Weil-Petersson metric along degenerations of Calabi-Yau manifolds}, Math. Res. Lett. {\bf 4} (1997), no. 1, 157--171.
\bibitem{We} J. Wehler, {\em Isomorphie von Familien kompakter komplexer Mannigfaltigkeiten}, Math. Ann. {\bf 231} (1977/78), no. 1, 77--90.
\bibitem{Ya}  S.-T. Yau, {\em On the Ricci curvature of a compact K\"ahler manifold and the complex Monge-Amp\`ere equation, I}, Comm. Pure Appl. Math. {\bf 31} (1978), 339--411.
\bibitem{Ya2}  S.-T. Yau, {\em A general Schwarz lemma for K\"ahler manifolds}, Amer. J. Math. {\bf 100} (1978), no. 1, 197--203.
\bibitem{ZT} Y. Zhang, {\em Convergence of K\"ahler manifolds and calibrated fibrations}, PhD thesis, Nankai Institute of Mathematics, 2006.
\bibitem{Zh} Y. Zhang, {\em Completion of the moduli space for polarized Calabi-Yau manifolds}, J. Differential Geom. {\bf 103} (2016), no. 3, 521--544.
\bibitem{Zh2} Y. Zhang, {\em Degeneration of Ricci-flat Calabi-Yau manifolds and its applications}, in {\em Uniformization, Riemann-Hilbert correspondence, Calabi-Yau manifolds \& Picard-Fuchs equations}, 551--592, Adv. Lect. Math. (ALM), 42, Int. Press, Somerville, MA, 2018.
\bibitem{Zh3} Y. Zhang, {\em Note on equivalences for degenerations of Calabi-Yau manifolds}, in {\em Surveys in Geometric Analysis 2017}, 186--202, Science Press, Beijing, 2018.
 \end{thebibliography}
\end{document}